\documentclass[10pt]{article}
\usepackage{amsmath}
\usepackage{amsthm}
\usepackage{amssymb}
\usepackage{amscd}
\usepackage{amsfonts}  
\usepackage{hyperref}
\usepackage[all]{xy}
\usepackage{graphicx}    
\usepackage{multicol}    

\DeclareMathAlphabet\EuR{U}{eur}{m}{n}
\SetMathAlphabet\EuR{bold}{U}{eur}{b}{n}

\newcommand{\calfin}{{\mathcal F}{\mathcal I}{\mathcal N}}
\newcommand{\calvcyc}{{\mathcal V}{\mathcal C}{\mathcal Y}{\mathcal C}}

\newcommand{\calall}{{\mathcal A}{\mathcal L}{\mathcal L}}

\newcommand{\calc}{{\cal C}}
\newcommand{\calcl}{{\cal C\!L}}
\newcommand{\cald}{{\cal D}}

\newcommand{\calf}{{\cal F}}
\newcommand{\calg}{{\cal G}}
\newcommand{\calh}{{\cal H}}

\newcommand{\cali}{{\cal I}}

\newcommand{\calp}{{\cal P}}

\newcommand{\calr}{{\cal R}}

\newcommand{\IK}{{\mathbb K}}

\newcommand{\IN}{{\mathbb N}}

\newcommand{\IQ}{{\mathbb Q}}
\newcommand{\IR}{{\mathbb R}}

\newcommand{\IZ}{{\mathbb Z}}

\newcommand{\bfE}{{\mathbf E}}
\newcommand{\bff}{{\mathbf f}}
\newcommand{\bfF}{{\mathbf F}}

\newcommand{\bfK}{{\mathbf K}}
\newcommand{\bfKH}{{\mathbf K}{\mathbf H}}
\newcommand{\bfL}{{\mathbf L}}

\newcommand{\bfNH}{{\mathbf N}{\mathbf H}}
\newcommand{\bfNil}{{\mathbf N}{\mathbf i}{\mathbf L}}

\newcommand{\bfT}{{\mathbf T}}

\newcommand{\bfid}{{\mathbf i}{\mathbf d}}

\newcommand{\curs}{\EuR}

\newcommand{\GROUPOIDS}{\curs{GROUPOIDS}}

\newcommand{\categoryN}{\curs{N}}
\newcommand{\NIL}{\curs{NIL}}
\newcommand{\Or}{\curs{Or}}

\newcommand{\SPACES}{\curs{SPACES}}
\newcommand{\SPECTRA}{\curs{SPECTRA}}

\newcommand{\asmb}{\operatorname{asmb}}

\newcommand{\colim}{\operatorname{colim}}

\newcommand{\fib}{\operatorname{fib}}

\newcommand{\id}{\operatorname{id}}
\newcommand{\ind}{\operatorname{ind}}

\newcommand{\Ker}{\operatorname{ker}}
\newcommand{\map}{\operatorname{map}}

\newcommand{\pr}{\operatorname{pr}}

\newcommand{\res}{\operatorname{res}}

\newcommand{\topo}{\operatorname{top}}

\newcommand{\Nil}{\operatorname{Nil}}  

\newcommand{\pt}{\operatorname{pt}}

\newcommand{\UNil}{\operatorname{UNil}}

\newcommand{\conering}{\Lambda}
\newcommand{\Map}{\operatorname{Map}}

\newcommand{\x}{\times}
\newcommand{\ox}{\otimes}

\newcommand{\KH}{\mathit{KH}}

\newcommand{\NK}{{\mathit{NK}}}

\newcommand{\hra}{\hookrightarrow}

\newcommand{\ULm}{{\underline{m}}}
\newcommand{\ULn}{{\underline{n}}}

\newcommand{\UL}{\underline}

\theoremstyle{plain}
\newtheorem{theorem}{Theorem}[section]
\newtheorem{lemma}[theorem]{Lemma}
\newtheorem{corollary}[theorem]{Corollary}
\newtheorem{proposition}[theorem]{Proposition}
\newtheorem{conjecture}[theorem]{Conjecture}

\theoremstyle{definition}
\newtheorem{definition}[theorem]{Definition}

\newtheorem{remark}[theorem]{Remark}

\theoremstyle{remark}

{\catcode`@=11\global\let\c@equation=\c@theorem}

\newcommand{\EGF}[2]{E_{#2}(#1)}               

\newcommand{\comsquare}[8]                   
{\begin{CD}
#1 @>#2>> #3\\
@V{#4}VV @VV{#5}V\\
#6 @>>#7> #8
\end{CD}
}

\newcommand{\comsquarel}[8]                   
{\begin{CD}
#1 @>#2>> #3\\
@V{#4}VV @V{#5}VV\\
#6 @>>#7> #8
\end{CD}
}

\newcommand{\xycomsquare}[8]                   
{\xymatrix
{
#1 \ar[r]^{#2} \ar[d]^{#4} &
#3 \ar[d]^{#5}  \\
#6\ar[r]^{#7} &
#8
}
}

\newcommand{\xycomsquarel}[8]                   
{\xymatrix
{
#1 \ar[r]^{#2} \ar[d]_{#4} &
#3 \ar[d]^{#5}  \\
#6\ar[r]^{#7} &
#8
}
}

\begin{document}

\title{Isomorphism Conjecture for homotopy $K$-theory
and groups acting on trees}
\author{
Arthur Bartels\thanks{bartelsa@math.uni-muenster.de}
and
Wolfgang L\"uck\thanks{\noindent email:
lueck@math.uni-muenster.de\protect\\
www: ~http://www.math.uni-muenster.de/u/lueck/\protect\\
FAX: 49 251 8338370\protect}
\\
Fachbereich Mathematik\\ Universit\"at M\"unster\\
Einsteinstr.~62\\ 48149 M\"unster\\Germany}

\maketitle

\typeout{-----------------------  Abstract  ------------------------}

\begin{abstract}
We discuss an analogon to the Farrell-Jones Conjecture for homotopy algebraic $K$-theory.
In particular, we prove that if a group $G$ acts on a tree and all isotropy groups satisfy this conjecture,
then $G$ satisfies this conjecture. This result can be used to get rational injectivity results for the
assembly map in the Farrell-Jones Conjecture in algebraic $K$-theory.
\smallskip


\end{abstract}


\typeout{------------ Introduction ---------------}
\setcounter{section}{-1}
\section{Introduction}
\label{sec: Introduction}

The Farrell-Jones Conjecture \cite{Farrell-Jones(1993a)}
in algebraic $K$-theory is concerned with the $K$-theory $K_n(RG)$ of  group rings
$RG$ for a group $G$ and a ring $R$.
The conjecture states that the assembly map
\begin{eqnarray}
H_n^G(\EGF{G}{\calvcyc};\bfK_R)  \to K_n(RG)
\label{assembly-map}
\end{eqnarray}
is an isomorphism. (This map is constructed
by applying a certain $G$-homology theory $H^G_n(-;\bfK_R)$ to the projection
$\EGF{G}{\calvcyc} \to \pt$, see
Definition~\ref{def: (Fibered) Isomorphism Conjectures for calh^?_*} and
Remark~\ref{rem: relation to the fibered isomorphism conjecture of Farrell-Jones}.)
There seem to occur two quite different phenomena in the algebraic $K$-theory
of such group rings. Firstly, $K_n(RG)$ contains elements coming from
the $K$-theory of $RF$ for finite subgroups $F$ of $G$. Secondly, it contains
nilgroup information. This is already illuminated in the simple case $G=\IZ$,
then $R[\IZ] = R[t,t^{-1}]$ and by the Bass-Heller-Swan
formula \cite{Bass-Heller-Swan(1964)}, \cite{Grayson(1976)}
\begin{eqnarray}
K_n(R[\IZ]) & \cong & K_n(R) \oplus K_{n-1}(R) \oplus \NK_{n}(R) \oplus \NK_{n}(R).
\label{Bass-Heller-Swan-formula}
\end{eqnarray}
Here $\NK_{n}(R)$ are the $\Nil$-groups of $R$, which can be defined
as the kernel of the projection
$K_n(R[t]) \to K_n(R)$ induced from $t \mapsto 0$.
In general it is known \cite{Bartels(2003b)}
that the domain of the assembly map \eqref{assembly-map} splits as
\begin{eqnarray}
H_n^G(\EGF{G}{\calfin};\bfK_R) \oplus
H_n^G(\EGF{G}{\calvcyc},\EGF{G}{\calfin} ;\bfK_R).
\end{eqnarray}
Thus, the Farrell-Jones Conjecture predicts a similar splitting for $K_n(RG)$.

In this paper we will formulate a (Fibered) Isomorphism Conjecture for homotopy algebraic $K$-theory,
see Conjecture~\ref{con: KH-Isomorphism-Conjecture}.
This variant of $K$-theory was defined by Weibel \cite{Weibel(1989)}, building on
the definition of Karoubi-Villamayor $K$-theory.
The homotopy algebraic $K$-theory groups of a ring $R$ are denoted by
$\KH_n(R)$. Their crucial property is homotopy invariance:
$\KH_n(R) \cong \KH_n(R[t])$. In particular, homotopy algebraic
$K$-theory does not contain $\Nil$-groups. We think about this
$\KH$-Isomorphism Conjecture as an
Isomorphism Conjecture for algebraic $K$-theory {\em modulo $\Nil$-groups}.
For a more precise formulation of the relation of the
Farrell-Jones Conjecture in algebraic $K$-theory to the
$\KH$-Isomorphism Conjecture see
Section~\ref{sec: The Relation between the K- and the KH-Isomorphism Conjecture}.

Our main results concerning the $\KH$-Isomorphism Conjecture are inheritance properties.
A group $G$ acts on a tree $T$, if $T$ is a $1$-dimensional
$G$-$CW$-complex which is contractible (after forgetting the group action).

\begin{definition}[The class of groups $\calc_0$]  \label{def: class C_0}
We define the following properties a class $\calc$ of groups
may or may not have.

\begin{description}

\item [(FIN)] All finite groups belong to $\calc$;

\item [(TREE)] Suppose that $G$ acts on a tree $T$.
Assume that for each $x \in T$ the isotropy
group $G_x$ belongs to $\calc$. Then $G$ belongs to $\calc$;

\item[(COL)]  Let $G$ be a group  with a directed system of subgroups
$\{G_i \mid i \in I\}$,
which is directed by inclusion and satisfies $\bigcup_{i \in I} G_i =
G$. If each $G_i$ belongs to $\calc$, then $G \in \calc$;

\item[(SUB)] If $G \in \calc$ and $H \subseteq G$ is a subgroup, then $H \in \calc$.

\end{description}
We define $\calc_0$ to be the smallest class of groups satisfying (FIN), (TREE) and (COL).
\end{definition}

It is not hard to check that the class $\calc_0$ is closed under taking subgroups.
For instance if $H$ is a subgroup of a group $G$ acting on a tree,
then $H$ acts also on this tree and the isotropy groups satisfy $H_x \subseteq G_x$.
By induction we may assume that the $G_x$ are closed under taking
subgroups and therefore $H \in \calc_0$.

\begin{theorem}{\bf (Inheritance properties of the $\KH$-Isomorphism Conjecture)}
\label{the: inhertitance for KH}
The class of groups  satisfying the Fibered $\KH$-Isomorphism
Conjecture for a fixed coefficient ring $R$ has the properties (FIN), (TREE), (COL) and (SUB).
The class of groups satisfying the $\KH$-Isomorphism
Conjecture for a fixed coefficient ring $R$ has the properties (FIN), (TREE) and (COL).
In particular, all groups in $\calc_0$ satisfy the (Fibered) $\KH$-Isomorphism Conjecture.
\end{theorem}

\begin{remark} \label{rem: extending calc}
The class of groups satisfying the $\KH$-Isomorphism Conjecture is stricly bigger than $\calc_0$ since it contains all
fundamental groups of closed Riemannian manifolds with negative sectional curvature
by \cite{Bartels-Reich(2005)} and
Theorem~\ref{the: K-theory version implies KH-version}~\ref{the: K-theory version implies KH-version (polynom)}.
\end{remark}

This result has the following applications.

\begin{theorem}[Extensions of groups and actions on trees]
\label{the: extensions and actions with finite stabilizers}
Let $1 \to K \to G \to Q \to 1$ be an extension of groups.
Suppose that $K$ acts on a tree with finite stabilizers
and that $Q$ satisfies the Fibered $\KH$-Isomorphism
Conjecture~\ref{con: KH-Isomorphism-Conjecture}  for the ring $R$. Then $G$ satisfies the Fibered
$\KH$-Isomorphism Conjecture~\ref{con: KH-Isomorphism-Conjecture} for the ring $R$.
\end{theorem}

A ring $R$ is called \emph{regular} if it is Noetherian and every
finitely generated $R$-module possesses a finite-dimensional
resolution by finitely generated projective modules.

\begin{theorem}{\bf (Conclusions for the $K$-theoretic Farrell-Jones Conjecture for groups in $\calc$).}
\label{the: conclusions for calc and K-theoretic Farrell-Jones}
Let $G$ be a group in the class $\calc_0$ defined above in \eqref{def: class C_0}. Then
\begin{enumerate}

\item Let $R$ be a regular ring with $\IQ \subseteq R$. Then the assembly map
$$H_n^G(\EGF{G}{\calfin};\bfK_R) \to K_n(RG)$$
is injective, or, equivalently, the injectivity part of the
Farrell-Jones Isomorphism Conjecture for algebraic $K$-theory is true
for $(G,R)$.

\item Let $R$ be the ring $\IZ$ of integers. Then the assembly map
$$H_n^G(\EGF{G}{\calfin};\bfK_{\IZ}) \to K_n(\IZ G)$$
is rationally injective,
or, equivalently, the rational injectivity part of the
Farrell-Jones Isomorphism Conjecture for algebraic $K$-theory is true for $(G,\IZ)$.
\end{enumerate}

\end{theorem}

\begin{proposition} \label{pro: groups belonging to calc}

The following classes of  groups belong to $\calc_0$:

\begin{enumerate}

\item \label{pro: groups belonging to calc: one relator groups} One relator groups;

\item \label{pro: groups belonging to calc: poly-free}
$G$ is poly-free, i.e. there is a filtration
$$\{1\} = G_0 \subseteq G_1 \subseteq G_2 \subseteq \ldots \subseteq G_n = G$$
such that $G_i$ is normal in $G_{i+1}$  with a free group as quotient $G_{i+1}/G_i$
The pure braid group is an example;

\item \label{pro: groups belonging to calc: 3-amifold groups}
Let $M$ be a compact orientable $3$-manifold with prime decomposition
$M = M_1 \sharp M_2 \sharp  \ldots \sharp M_n$. Suppose that
each $M_i$, which has infinite fundamental group and is aspherical, has
a boundary or is a Haken manifold. Then $\pi_1(M) \in \calc_0$;

\item \label{pro: groups belonging to calc: surface groups}
If $M$ is a compact $2$-dimensional manifold, then $\pi_1(M) \in \calc_0$;

\item \label{pro: groups belonging to calc: submanifolds of S^3}
If $M$ is a submanifold of $S^3$, then $\pi_1(M) \in \calc_0$.

\end{enumerate}

\end{proposition}

Next we discuss similar inheritance properties
for the Farrell-Jones Conjecture in algebraic K-theory.
A ring $R$ is called \emph{regular coherent} if every finitely presented
$R$-module possesses a finite-dimensional resolution by finitely
generated projective $R$-modules.
A ring $R$ is regular if and only if it is regular coherent and Noetherian.
A group $G$ is called \emph{regular} or \emph{regular coherent} respectively
if for any regular ring $R$ the group ring $RG$ is regular respectively regular coherent.
For more information about these notions we refer to \cite[Theorem 19.1]{Waldhausen(1978a)}.

\begin{definition}[The classes of groups $\calcl$ and $\calcl'$] \label{def: class CL}
Consider the following further properties a class $\calc$ of groups may have.
\begin{description}
\item[(TRI)]  The trivial group belongs to  $\calc$;
\item[(VCYC)] All virtually cyclic groups belong to $\calc$;
\item[(TREE$_\calr$)] Suppose that $G$ acts on a tree $T$.
              Assume that for each $x \in T$ the isotropy
              group $G_x$ belongs to $\calc$. For each edge $e$ of $T$, assume that
              the isotropy group $G_e$ is regular coherent.
              Then $G$ belongs to $\calc$;
\end{description}
The class $\calcl$ is defined as the smallest class of groups
satisfying (TRI), (TREE$_\calr$) and (COL).
The class $\calcl'$ is defined as the smallest class of groups
satisfying (VCYC), (TREE$_\calr$) and (COL).
\end{definition}

All groups appearing in $\calcl$ are torsionfree.
Similar to the class $\calc_0$,
the classes $\calcl$ and $\calcl'$ are closed under taking subgroups \cite[Proposition 19.3]{Waldhausen(1978a)}.
We conclude from Waldhausen~\cite[Theorem 17.5 on page 250]{Waldhausen(1978a)}
that $\calcl$ contains a group $G$  appearing in Proposition~\ref{pro: groups belonging to calc}
under \ref{pro: groups belonging to calc: one relator groups}, \ref{pro: groups belonging to calc: 3-amifold groups},
\ref{pro: groups belonging to calc: surface groups}, and \ref{pro: groups belonging to calc: submanifolds of S^3}
provided that $G$ is torsionfree.
One of the main results in Waldhausen's article  \cite{Waldhausen(1978a)} is that for a regular ring $R$
the $K$-theoretic assembly map
$$
H_n(BG;\bfK_R) \to K_n(RG)
$$
is an isomorphism.
Actually, Waldhausen states this only for $n \geq 0$, but the embedding of $K_{n-1}(R)$ into $K_n(R[\IZ])$
allows the extension to all $n$, see for example Remark~\ref{rem: nil vanishes(amalgamation)}.
Furthermore, Waldhausen considers HNN-extensions and amalgamated products rather than action on trees, but
this does not change the class $\calcl$, compare Remark~\ref{rem: properties (TREE) and (TREE')} and
Lemma~\ref{lem: tree property for calh^?_*(-,bfE)}.

\begin{theorem}
\label{the: inhertiance for Farrell Jones}
Let $R$ be a regular ring.
The class of groups satisfying the Farrell-Jones Conjecture in algebraic
$K$-theory for the ring $R$ has the properties (VCYC), (TREE$_\calr$) and (COL).
In particular, all groups in $\calcl'$ satisfy the Farrell-Jones Conjecture in algebraic $K$-theory for
the ring $R$.
\end{theorem}

Related results can be found in \cite{Juan-Pineda-Prassidis(2004)} and \cite{Roushon(2004)}.
It is an interesting question, whether the class of groups satisfying the Farrell-Jones Conjecture
in algebraic $K$-theory has the property (TREE).

Theorem~\ref{the: inhertitance for KH} and
Waldhausen's result imply

\begin{corollary} \label{K cong KH for g in CL}
Let $R$ be a regular ring. Let $G$ be a group in the class $\calcl$. Then
the canonical  map
$$K_n(RG) \to \KH_n(RG)$$
is bijective for $n \in \IZ$.
\end{corollary}

Isomorphism Conjectures can be formulated in the quite general context
of equivariant homology theories, see
Definition~\ref{def: (Fibered) Isomorphism Conjectures for calh^?_*}.
We show in Theorem~\ref{the: Isomorphisms Conjecture and actions on trees}
that the property (TREE) holds for the class of groups satisfying the
Isomorphism Conjecture for such an equivariant homology theory
whenever the equivariant homology theory satisfies
the \emph{tree property}, see Definitions~\ref{def: tree-property}.
The weaker property (TREE$_\calr$) is related to the \emph{regular} tree property
(see Definition~\ref{def: tree-property}), which is a weakening of the tree property.
The above Theorem~\ref{the: extensions and actions with finite stabilizers}
has also an analogon in this setting, see
Corollary~\ref{cor: extensions by groups acting on trees with finite stabilizers}.
The tree property means essentially that there are Mayer-Vietoris
sequences for amalgamated products and
HNN extensions  of groups in the equivariant homology theory
(see  Lemma~\ref{lem: suffices to check tree property on amal. prod. and HNN}).
On spectrum level this means that there are certain homotopy cartesian diagrams of spectra
(see Lemma~\ref{lem: tree property for calh^?_*(-,bfE)} and Remark~\ref{rem: smaller diagrams}).
In Section~\ref{sec: The KH-Isomorphism Conjecture}
we define the equivariant homology theory $H_*^?(-,\bfKH_R)$ that is
relevant for the $\KH$-Isomorphism Conjecture.
We prove in Theorem~\ref{the: tree-property-for-KH}
that this theory satisfies the tree property.
In the case of algebraic $K$-theory amalgamated free products and
HNN extensions have been analyzed by Waldhausen \cite{Waldhausen(1978a)}.
In both cases there are long exact sequences, but they involve as an additional term Waldhausen's $\Nil$-groups.
Their nontriviality obstructs the equivariant homology theory $H^?_*(-;\bfK_R)$ relevant for the Farrell-Jones
Conjecture in algebraic $K$-theory from having the tree property.
Our proof of the tree property for $H_*^?(-,\bfKH_R)$, consists essentially of showing that Waldhausen's $\Nil$-groups
are killed under the transition from $K$ to $\KH$.
On the other hand vanishing results for Waldhausen's $\Nil$-groups can be used to show
(see Theorem~\ref{the: regular tree property for K})
that $H^?_*(-;\bfK_R)$ has the \emph{regular} tree property (see Definition~\ref{def: tree-property}).
This implies then, that for a regular ring $R$ the class of groups satisfying the Farrell-Jones Conjecture has
the property (TREE$_\calr$), see
Theorem~\ref{the: Isomorphisms Conjecture and actions on trees}~
\ref{the: Isomorphisms Conjecture and actions on trees: vcyc}.
It is an interesting question for which rings $R$ the equivariant
homology theory $H_*^?(-;\bfK_R) \ox \IQ$ has the tree property.
It is worthwhile to consider also $H_*^?(-;\bfL^{-\infty}_R)$,
the equivariant homology theory
relevant for the Farrell-Jones Conjecture in $L$-theory.
In this case amalgamated free products
and HNN extensions have been analyzed by Cappell~\cite{Cappell(1974c)}.
Again additional terms appear in the long exact sequences,
the $\UNil$-groups and non-triviality of those obstructs this theory from having the tree property.
On the other hand these $\UNil$-groups are known to be
$2$-torsion \cite{Cappell(1974c)},
thus $H_n^?(-;\bfL^{-\infty}_R) \ox \IZ[\frac{1}{2}]$ does have the tree
property. Thus we obtain the following result.

\begin{theorem}{\bf (Conclusions for the $L$-theoretic Farrell-Jones Conjecture for groups in $\calc$).}
 \label{the: conclusions for calc and L-theoretic Farrell-Jones}
The class of groups for which the assembly map
\[
H^G_n(\EGF{G}{\calfin};\bfL^{-\infty}_R) \to L^{-\infty}_n(RG)
\]
becomes an isomorphism after tensoring with $\IZ[\frac{1}{2}]$,
has the properties (FIN), (TREE) and (COL).
In particular, this class contains the class in $\calc_0$ from Definition~\ref{def: class C_0}.
\end{theorem}

In the context of topological
$K$-theory, i.e. for
the Baum-Connes Conjecture one can apply our results to the
equivariant $K$-theory $H^G_n(-;\bfK^{\topo}) = K^G_n(-)$.
Then one obtains the analogon of our Theorem~\ref{the: inhertitance for KH}.
In this case amalgamated products and HNN extensions have been
analyzed Pimsner-Voiculescu \cite{Pimsner-Voiculescu(1982)} and Pimsner
\cite{Pimsner(1986)}. Here the situation is much better, since no $\Nil$-groups appear.
This analogon has already been proved by Oyono-Oyono \cite{Oyono-Oyono(2001b)}
for the Baum-Connes Conjecture (with coefficients).

We are indebted to Holger Reich for pointing out the reference
\cite{Weibel(1981)} to us.

The papers is organized as follows:

\begin{tabular}{lll}
\ref{sec: Isomorphism Conjectures for equivariant homology theories}. & &
Isomorphism Conjectures for equivariant homology theories
\\
\ref{sec: Homological aspects}. & & Homological aspects
\\
\ref{sec: Continuous equivariant homology theories} & &
Continuous equivariant homology theories
\\
\ref{sec: The tree property}. & &
The tree property
\\
\ref{sec: Equivariant homology theories constructed from spectra}.
& &
Equivariant homology theories constructed from spectra
\\
\ref{sec: Isomorphism Conjectures for spectra}. & &
Isomorphism Conjectures for spectra
\\
\ref{sec: The KH-Isomorphism Conjecture}. & &
The $\KH$-Isomorphism Conjecture
\\
\ref{sec: The Relation between the K- and the KH-Isomorphism Conjecture}.
& &
The Relation between the $K$- and the $\KH$-Isomorphism Conjecture
\\
\ref{sec: Non-connective Waldhausen Nil} & &
Non-connective Waldhausen Nil
\\
\ref{subsec: Waldhausen's cartesian squares}. & &
Waldhausen's cartesian squares
\\
\ref{sec: The tree property for KH}. & &
The tree property for KH
\\
& & References
\end{tabular}


\typeout{----------- Isomorphisms Conjecture for equivariant homology theories   ----------------}
\section{Isomorphism Conjectures for equivariant homology theories}

\label{sec: Isomorphism Conjectures for equivariant homology theories}

We will use the notion of an equivariant homology theory $\calh^?_*$ with values in
$\Lambda$-modules for a commutative associative ring $\Lambda$ with unit from
\cite[Section 1]{Lueck(2002b)}.
This essentially means that we get for each group $G$ a $G$-homology theory $\calh^G_*$ which
assigns to a (not necessarily proper or cocompact) pair of
$G$-$CW$-complexes $(X,A)$ a $\IZ$-graded $\Lambda$-module $\calh^G_n(X,A)$, such that there exists natural
long exact sequences of pairs and $G$-homotopy invariance, excision, and the disjoint union axiom are satisfied.
Moreover, an induction structure is required which in particular implies for a subgroup
$H \subseteq G$ and a $H$-$CW$-pair $(X,A)$ that there is a natural isomorphism
$\calh^H_n(X,A) \xrightarrow{\cong} \calh^G_n(G \times_H(X,A))$.

We  will later discuss examples, the most important ones will be given by those equivariant
homology theories which appear in the Baum-Connes Conjecture and the Farrell-Jones Conjecture.
These conjectures are special cases of the following more
general formulation of a (Fibered) Isomorphism Conjecture
(see Section~\ref{sec: Isomorphism Conjectures for spectra}).

A family $\calf$ of subgroups of $G$ is a set of subgroups which is
closed under conjugation and taking subgroups.
If $\calc$ is a class of groups that is closed under taking subgroups
and isomorphisms, then the collections of subgroups of $G$ that are in $\calc$
forms a family $\calc(G)$ of subgroups of $G$.
Abusing notation, we will denote this family often by $\calc$.
Examples are the families $\calfin$ of finite subgroups
and $\calvcyc$ of virtually cyclic subgroups.
Given a group homomorphism
$\phi \colon K \to G$ and a family $\calf$ of subgroups of $G$,
define the family
$\phi^*\calf$ of subgroups of $K$ by
$\phi^*\calf = \{H \subseteq K \mid \phi(H) \in \calf\}$.
If $i\colon H \to G$ is the inclusion of a subgroup, then we write often $\calf|_H$ for $i^*\calf$.
Associated to such a family there is a $G$-$CW$-complex
$\EGF{G}{\calf}$ (unique up to $G$-homotopy equivalence)
with the property that the fixpoint sets $\EGF{G}{\calf}^H$
are contractible for $H \in \calf$ and
empty for $H \notin \calf$.
For $\calf = \calall$ the family of all subgroups,
we can take the one-point-space $\pt$
as a model for $\EGF{G}{\calall}$.
For more information about the spaces we refer for instance to \cite{Lueck(2004a)}.

\begin{definition}[(Fibered) Isomorphism Conjecture for $\calh^?_*$]
\label{def: (Fibered) Isomorphism Conjectures for calh^?_*}
Let $\calh^?_*$ be an equivariant homology theory with values in
$\Lambda$-modules. A group $G$ together with a family
of subgroups $\calf$ satisfies the
\emph{Isomorphism Conjecture (in the range $\le N$)}
if the projection $\pr \colon \EGF{G}{\calf} \to \pt$
to the one-point-space $\pt$
induces an isomorphism
$$\calh^G_n(\pr) \colon \calh^G_n(\EGF{G}{\calf})
     \xrightarrow{\cong} \calh^G_n(\pt)$$
for $n \in \IZ$  (with $n \le N)$.

The pair $(G,\calf)$  satisfies  the \emph{Fibered Isomorphism
Conjecture (in the range $\le N$)} if for each group homomorphism
$\phi \colon K \to G$ the pair $(K,\phi^*\calf)$ satisfies the
Isomorphism Conjecture (in the range $\le N$).
\end{definition}



Built in into the Fibered Isomorphism Conjecture is the following obvious inheritance
property which is not known to be true in general in the non-fibered case.

\begin{lemma} \label{lem: basic inheritance property of fibered conjecture}
Let $\calh^?_*$ be an equivariant homology theory,
let $\phi \colon K \to G$ be a group
homomorphism and let $\calf$ be a family of subgroups.
If $(G,\calf)$ satisfies
the Fibered Isomorphism Conjecture
\ref{def: (Fibered) Isomorphism Conjectures for calh^?_*} (in the range $\le N$),
then $(K,\phi^*\calf)$ satisfies
the Fibered Isomorphism Conjecture
\ref{def: (Fibered) Isomorphism Conjectures for calh^?_*} (in the range $\le N$).
\end{lemma}
\begin{proof}
If $\psi \colon L \to K$ is a group homomorphism, then
$\psi^*(\phi^*\calf) = (\phi \circ \psi)^*\calf$.
\end{proof}

In particular, if for a given class of groups $\calc$, which is closed under isomorphism and taking subgroups,
the Fibered Isomorphism Conjecture~\ref{def: (Fibered) Isomorphism Conjectures for calh^?_*}
is true for $(G,\calc(G))$ and if $H \subseteq G$ is a subgroup, then
the Fibered Isomorphism Conjecture~\ref{def: (Fibered) Isomorphism Conjectures for calh^?_*}
is true for $(H,\calc(H))$.


\typeout{------------- Homological aspects ----------------}
\section{Homological aspects}
\label{sec: Homological aspects}

The disjoint union axiom ensures that a $G$-homology is compatible with directed colimits.

\begin{lemma} \label{lem: G-homology theory and colimit}
Let $\calh^G_*$ be a $G$-homology theory.
Let $X$ be a $G$-$CW$-complex and $\{X_i \mid i \in I\}$ be a directed system of
$G$-$CW$-subcomplexes directed by inclusion such that $X = \bigcup_{i \in I} X_i$.
Then for all $n \in \IZ$ the natural map
$$\colim_{i \in I} \calh_n^G(X_i) \xrightarrow{\cong}   \calh^G_n(X)$$
is bijective.
\end{lemma}

\begin{proof} Compare for example with \cite[Proposition 7.53 on page 121]{Switzer(1975)},
where the non-equivariant case for $I=\IN$ is treated. The main point is
that the functor colimit over a directed system of $R$-modules is exact.
\end{proof}

\begin{lemma} \label{lem: EGF{G}{calf} times Z to Z}
Let $\calh^?_*$ be an equivariant homology theory with values in
$\Lambda $-modules in the sense of \cite[Section 1]{Lueck(2002b)}.
Let $G$ be a group and let $\calf$ be a family of subgroups of $G$.
Let $Z$ be a $G$-$CW$-complex. Consider $N \in \IZ \cup \{\infty\}$.
Suppose for each $H \subseteq G$ which occurs as isotropy group in $Z$
that the $G$-map induced by the projection
$\pr\colon \EGF{H}{\calf|_H} \to \pt$
$$\calh^H_n(\pr) \colon \calh^H_n(\EGF{H}{\calf|_H}) \to \calh^H_n(\pt)$$
is bijective for all $n \in \IZ, n \le N$.

Then the map induced by the projection
$\pr_1 \colon \EGF{G}{\calf} \times Z \to Z$
\begin{eqnarray*} \calh^G_n(\pr_1) \colon \calh^G_n(\EGF{G}{\calf} \times Z)
& \to & \calh^G_n(Z)
\end{eqnarray*}
is bijective for $n \in \IZ, n \le N$.
\end{lemma}
\begin{proof}
We first prove the claim for
finite-dimensional $G$-$CW$-complexes by induction over $d = \dim(Z)$.
The induction beginning $\dim(Z) = -1$, i.e. $Z = \emptyset$, is trivial.
In the induction step from $(d-1)$ to $d$ we choose a $G$-pushout
$$\comsquare{\coprod_{i \in I_d} G/H_i \times S^{d-1}}{}{Z_{d-1}}
{}{}
{\coprod_{i \in I_d} G/H_i \times D^d}{}{Z_d}
$$
If we cross it with $\EGF{G}{\calf}$, we obtain another
$G$-pushout of $G$-$CW$-complexes.
The various projections induce a map
from the Mayer-Vietoris sequence of the latter $G$-pushout to the
Mayer-Vietoris sequence of the
first $G$-pushout. By the Five-Lemma it suffices to prove
that the following maps
\begin{eqnarray*}
\calh_n^G(\pr_2) \colon \calh^G_n\left(\EGF{G}{\calf}
  \times \coprod_{i \in I_d} G/H_i \times S^{d-1}\right)
& \to &
\calh^G_n\left(\coprod_{i \in I_d} G/H_i \times S^{d-1}\right);
\\
\calh_n^G(\pr_3) \colon \calh^G_n(\EGF{G}{\calf} \times Z_{d-1})
& \to &  \calh^G_n(Z_{d-1});
\\
\calh_n^G(\pr_4) \colon \calh^G_n\left(\EGF{G}{\calf}
 \times \coprod_{i \in I_d} G/H_i \times D^n\right) & \to &
\calh^G_n\left(\coprod_{i \in I_d}G/H_i \times D^n\right),
\end{eqnarray*}
are bijective for $n \in \IZ, n \le N$. This follows from
the induction hypothesis for the first two maps.
Because of the disjoint union axiom  and $G$-homotopy invariance
of $\calh^?_*$ the claim follows for the third map
if we can show for any $H\subseteq G$ which occurs as isotropy group in
$Z$ that the map
\begin{eqnarray}
\calh_n^G(\pr_1) \colon \calh^G_n(\EGF{G}{\calf} \times G/H )
& \to &  \calh_n^G(G/H)
\label{calh^G_n(pr_1) for G/H}
\end{eqnarray}
is bijective for $n \in \IZ, n \le N$. The $G$-map
$$G \times_H \res_G^H \EGF{G}{\calf} \to G/H \times \EGF{G}{\calf}
  \hspace{5mm} (g,x) ~ \mapsto ~ (gH,gx)$$
is a $G$-homeomorphism where $\res_G^H$ denotes the restriction of
the $G$-action to an $H$-action.
Obviously $\res_G^H \EGF{G}{\calf}$ is a model for $\EGF{H}{\calf|_H}$.
Since for any $H$-$CW$-complex $Y$ there is a  natural isomorphism
$\calh^H_n(Y) \xrightarrow{\cong} \calh^G_n(G \times_HY)$,
the map \eqref{calh^G_n(pr_1) for G/H}  can be identified with the map
\begin{eqnarray*}
\calh_n^G(\pr) \colon \calh^H_n(\EGF{H}{\calf|_H}) & \to &  \calh^H_n(\pt)
\end{eqnarray*}
which is bijective for all $n \in \IZ, n \le N$ by assumption.
This finishes the proof
in the case  that $Z$ is finite-dimensional.

Finally we consider an arbitrary $G$-$CW$-complex $Z$.
It can be written as the colimit $\colim_{d \to \infty} Z_d$.
The natural maps
\begin{eqnarray*}
\colim_{d \to \infty} \calh^G_n(\EGF{G}{\calf} \times Z_d) &
\xrightarrow{\cong} & \calh^G_n(\EGF{G}{\calf} \times Z);
\\
\colim_{d \to \infty} \calh^G_n(Z_d) &  \xrightarrow{\cong} & \calh^G_n(Z),
\end{eqnarray*}
are bijective by Lemma~\ref{lem: G-homology theory and colimit}.
Since the colimit of isomorphisms
is an isomorphism again, Lemma \ref{lem: EGF{G}{calf} times Z to Z} follows.
\end{proof}

\begin{theorem}[Reducing the family]
\label{the: Reducing the family}
Let $\calh^?_*$ be an equivariant homology theory with values in $\Lambda $-modules.
Let $G$ be a group and let $\calf \subseteq \calg$ be families of subgroups of $G$.
Consider $N \in \IZ \cup \{\infty\}$.
Suppose for each $H \in \calg$, or, more generally, suppose for each isotropy group appearing in a
specific model for $\EGF{G}{\calg}$ that $(H,\calf|_H)$ satisfies the
(Fibered) Isomorphism Conjecture~\ref{def: (Fibered) Isomorphism Conjectures for calh^?_*} (in the range $\leq N$).

Then $(G,\calg)$ satisfies the
(Fibered) Isomorphism Conjecture~\ref{def: (Fibered) Isomorphism Conjectures for calh^?_*} (in the range $\leq N$)
if and only if $(G,\calf)$ satisfies the
(Fibered) Isomorphism Conjecture~\ref{def: (Fibered) Isomorphism Conjectures for calh^?_*} (in the range $\leq N$).
\end{theorem}
\begin{proof}
For the Isomorphism Conjecture
this follows from Lemma \ref{lem: EGF{G}{calf} times Z to Z}
applied to the case
$Z = \EGF{G}{\calg}$ and the fact that
$\EGF{G}{\calf} \times \EGF{G}{\calg}$ is a model for
$\EGF{G}{\calf}$.
The case of the Fibered Isomorphism Conjecture is easily reduced to the former case.
\end{proof}

\begin{lemma} \label{lem: fibered conjecture and extensions}
Let $\calh^?_*$ be an equivariant homology theory with values in $\Lambda $-modules.
Let $\calc$ be a class of groups that is closed under isomorphisms, subgroups and quotients.
Let $1 \to L \to G \xrightarrow{p} Q \to 1$ be an extension of groups.
Suppose that $(Q;\calc(Q))$ satisfies the
Fibered Isomorphism Conjecture \ref{def: (Fibered) Isomorphism Conjectures for calh^?_*} (in the range $\le N$)
and that for $H \in  p^*\calc(Q)$ the pair
$(H,\calc(H))$ satisfies the
Fibered Isomorphism Conjecture \ref{def: (Fibered) Isomorphism Conjectures for calh^?_*} (in the range $\le N$).

Then $(G,\calc(G))$ satisfies the
Fibered Isomorphism Conjecture \ref{def: (Fibered) Isomorphism Conjectures for calh^?_*} (in the range $\le N$).
\end{lemma}
\begin{proof}
By Lemma~\ref{lem: basic inheritance property of fibered conjecture} the pair $(G, p^*\calc(Q))$ satisfies the
Fibered Isomorphism Conjecture \ref{def: (Fibered) Isomorphism Conjectures for calh^?_*} (in the range $\le N$).
Since $\calc$ is closed under quotients we have $\calc(G) \subseteq p^*\calc(Q)$.
Now the assumption on the subgroups $H \in p^*\calc(Q)$ and Theorem~\ref{the: Reducing the family} imply the result.
\end{proof}


\typeout{------- Continuous equivariant homology theories  ---------------}
\section{Continuous equivariant homology theories}
\label{sec: Continuous equivariant homology theories}

In this section  we explain a criterion for an equivariant homology theory
ensuring that for a class  $\calc$ of groups closed under subgroups and isomorphisms
the (Fibered) Isomorphism Conjecture~\ref{def: (Fibered) Isomorphism Conjectures for calh^?_*} is true
for $(G,\calc(G))$ provided that $G$ is a directed union $G = \bigcup_{i \in I} G_i$ of groups $G_i$
and the (Fibered) $\calc$-Isomorphism Conjecture~\ref{def: (Fibered) Isomorphism Conjectures for calh^?_*} is true
for $(G_i,\calc(G_i))$ for all $i \in I$.

\begin{definition}[Continuous equivariant homology theory]
\label{def: continuous equivariant homology theory}
An equivariant homology theory $\calh^?_*$ is called \emph{continuous}
if for each group $G$  and directed system of subgroups
$\{G_i \mid i \in I\}$,
which is directed by inclusion and satisfies
$\bigcup_{i \in I} G_i = G$, and each  $n \in \IZ$ the map
$$\colim_{i \in I} j_i \colon \colim_{i \in I} \calh^{G_i}_n(\pt)
      \to \calh^G_n(\pt)$$
is an isomorphism, where
$j_i \colon \calh^{G_i}_n(\pt) \to \calh^G_n(\pt)$ is the composition
of the induction isomorphism
$\calh^{G_i}_n(\pt) \xrightarrow{\cong}  \calh_n^G(G/G_i)$ with the map induced
by the projection $G/G_i \to \pt$.
\end{definition}

\begin{lemma}
\label{lem: colim G_i isomorphism for X}
Let $\calh_*^?$ be a continuous equivariant homology theory.
Let $G$ be a group with a directed system of
subgroups $\{G_i \mid i \in I\}$,
which is directed by inclusion and satisfies $\bigcup_{i \in I} G_i = G$.

Then for each $G$-$CW$-complex $X$ and each  $n \in \IZ$ the map
$$\colim_{i \in I} j_i \colon \colim_{i \in I} \calh^{G_i}_n(\res_G^{G_i} X)
        \to \calh^G_n(X)$$
is an isomorphism, where
$j_i \colon \calh^{G_i}_n(\res_G^{G_i} X ) \to \calh^G_n(X)$ is the composition
of the induction isomorphism $\calh^{G_i}_n(\res_G^{G_i} X)
\xrightarrow{\cong}  \calh_n^G(G \x_{G_i} \res_G^{G_i} X)$
with the homomorphism induced
by the $G$-map $G \x_{G_i} \res_G^{G_i} X \to X$ that sends $(g,x)$ to
$gx$.
\end{lemma}
\begin{proof}
Since $\colim_{i \in I}$ is an exact functor,
$\colim_{i \in I} \calh^{G_i}_n(\res_G^{G_i} X)$ is a
$G$-homology theory in $X$.
The map $\colim_{i \in I} j_i$ is a transformation of $G$-homology theories.
Therefore it suffices  to prove that
\begin{eqnarray*}
\colim_{i \in I} j_i \colon \colim_{i \in I} \calh^{G_i}_n(\res_G^{G_i} G/H) & \to & \calh^G_n(G/H)
\end{eqnarray*}
is an isomorphism for every subgroup $H \subseteq G$ and $n \in \IZ$.

For $i \in I$ let $k_i \colon G_i/G_i \cap H \to \res_G^{G_i} G/G_i$ be the obvious injective $G_i$ map.
Then the following diagram commutes
$$\xymatrix
{
\colim_{i \in I} \calh_n^{G_i \cap H}(\pt) \ar[r]^{\cong} \ar[dd]_{\colim_{i \in I} j_i}^{\cong} &
\colim_{i \in I} \calh_n^{G_i}(G_i/G_i \cap H) \ar[d]^{\colim_{i \in I} k_i} \\
& \colim_{i \in I} \calh_n^{G_i}(\res_G^{G_i}G/H) \ar[d]^{\colim_{i \in I} j_i} \\
\calh_n^H(\pt)  \ar[r]^{\cong} & \calh_n^G(G/H)
}
$$
where the horizontal maps are the isomorphism given by induction.
The left vertical arrow is bijective since
$\calh^?_*$ is continuous by assumption.
Hence it remains to show that the map
\begin{eqnarray}
\colim_{i \in I} k_i \colon \colim_{i \in I} \calh_n^{G_i}(G_i/G_i \cap H)
\to
 \colim_{i \in I} \calh_n^{G_i}(\res_G^{G_i}G/H) &&
\label{colim_{i in I} k_i}
\end{eqnarray}
is surjective.

We get an obvious decomposition of $G_i$-sets
$$\res_G^{G_i} G/H ~ = ~ \coprod_{G_igH \in G_i\backslash G/H}~ G_i/G_i \cap gHg^{-1}.$$
It induces an identification
$$\calh^{G_i}_n(\res_G^{G_i} G/H) ~ = ~ \bigoplus_{G_igH \in G_i\backslash G/H}~ \calh^{G_i}_n(G_i/G_i \cap gHg^{-1}).$$
The summand corresponding to $G_i1H$ is precisely the image of
$$\calh^{G_i}_n(k_i) \colon \calh^{G_i}_n(G_i/G_i \cap H) \to \calh^{G_i}_n(\res_G^{G_i} G/H).$$
Consider an element $G_igH \in G_i \backslash G/H$. Choose an index
$j$ with $j \ge i$ and $g \in G_j$. Then the structure map for $i \le j$
is a map $\calh^{G_i}_n(\res_G^{G_i} G/H) \to \calh^{G_j}_n(\res_G^{G_j} G/H)$ which sends the summand
corresponding to  $G_igH \in G_i \backslash G/H$ to the image of
$$\calh^{G_j}_n(k_j) \colon \calh^{G_j}_n(G_j/G_j \cap H) \to \calh^{G_j}_n(\res_G^{G_j} G/H).$$
This implies that the map \eqref{colim_{i in I} k_i}
is surjective. This finishes the proof of Lemma~\ref{lem: colim G_i isomorphism for X}.
\end{proof}

\begin{proposition}
\label{prop: isomorphism conjecture is stable under colim}
Let $\calh_*^?$ be an equivariant homology theory which is continuous.
Let $\calc$ be a class of groups that is closed under isomorphism and taking subgroups.
Let $G$ be the directed union $G = \bigcup_{i \in I} G_i$ of subgroups $G_i$
such that the (Fibered) Isomorphism Conjecture~\ref{def: (Fibered) Isomorphism Conjectures for calh^?_*}
(in the range $\leq N$)  is true
for $(G_i,\calc(G_i))$ for all $i \in I$.

Then  the (Fibered) Isomorphism Conjecture~\ref{def: (Fibered) Isomorphism Conjectures for calh^?_*}
(in the range $\leq N$)  is true for $(G,\calc(G))$.
\end{proposition}
\begin{proof}
Since $\calh^?_*$ is continuous by assumption, we get the isomorphism
$$\colim_{i \in I} \calh^{G_i}_*(\pt) = \calh_*^G(\pt)$$
and from Lemma~\ref{lem: colim G_i isomorphism for X} the isomorphism
$$\calh_*^G(\EGF{G}{\calc(G)}) = \colim_{i \in I} \calh_*^{G_i}(\res_G^{G_i} \EGF{G}{\calc(G)}).$$
The result follows for the Isomorphism Conjecture since the colimit of an
isomorphism is an isomorphism and since
$\res_G^{G_i} \EGF{G}{\calc(G)}$ is a model for $\EGF{G_i}{\calc(G_i)}$.
If $\phi\colon K \to G$ is a group homomorphism then the same argument can be applied to the triple
$(K,\phi^*\calc(G), \{ \phi^{-1}(G_i) \mid i \in I \} )$ in place of
$(G,\calc(G),\{G_i \mid i \in I\})$ and this
implies the statement for the Fibered Isomorphism Conjecture.
\end{proof}


\typeout{--------- The tree property  ----------------}
\section{The tree property}
\label{sec: The tree property}

In this section  we study criteria for an equivariant homology theory that ensure
that the class of groups $G$ for which $(G,\calfin)$ satisfies the (Fibered)
Isomorphism Conjecture~\ref{def: (Fibered) Isomorphism Conjectures for calh^?_*}
has property (TREE) from Definition~\ref{def: class C_0}
or that the class of groups $G$ for which $(G,\calvcyc)$ satisfies the
Isomorphism Conjecture~\ref{def: (Fibered) Isomorphism Conjectures for calh^?_*}
has property (TREE$_\calr$) from Definition~\ref{def: class CL}.

\begin{definition}[Tree Property] \label{def: tree-property}
An equivariant homology theory $\calh^?_*$ has the
\emph{tree property}
if for any group $G$ that acts on a tree $T$,
the projection $\pr \colon T \to \pt$ induces for all $n \in
\IZ$ (with $n \le N$) isomorphisms
$$\calh_n^G(\pr) \colon \calh^G_n(T) ~ \xrightarrow ~ \calh_n^G(\pt).$$
It has the
\emph{regular tree property}
if for any group $G$ that acts on a tree $T$, such that
for each edge $e$ of $T$ the isotropy group $G_e$ is regular coherent,
the projection $\pr \colon T \to \pt$ induces for all $n \in
\IZ$ (with $n \le N$) isomorphisms
$$\calh_n^G(\pr) \colon \calh^G_n(T) ~ \xrightarrow ~ \calh_n^G(\pt).$$
\end{definition}

\
\begin{theorem}
{\bf (The tree property and inheritance properties of Isomorphism Conjectures)}
\label{the: Isomorphisms Conjecture and actions on trees}
Let $\calh^?_*$ be an equivariant homology theory.
Let $\calc$ be a class of groups closed under subgroups and isomorphisms.
Let $\cald_{\fib}(\calc)$ be the class of groups $G$
for which Fibered Isomorphism Conjecture~\ref{def: (Fibered) Isomorphism Conjectures for calh^?_*}
(in the range $\leq N$) is true for $(G,\calc(G))$ and let $\cald(\calc)$ be the class of groups $G$
for which Isomorphism Conjecture~\ref{def: (Fibered) Isomorphism Conjectures for calh^?_*}
(in the range $\leq N$) is true for $(G,\calc(G))$.
\begin{enumerate}

\item \label{the: Isomorphisms Conjecture and actions on trees: fibered}
      Suppose that $\calh^?_*$ has the tree property \eqref{def: tree-property}. Then the class $\cald_ {\fib}(\calc)$ has
      the property (TREE) from Definition~\ref{def: class C_0}.

\item \label{the: Isomorphisms Conjecture and actions on trees: fin}
      Suppose $\calh^?_*$ has the tree property \eqref{def: tree-property} and $\calc \subseteq \calfin$.
      Then the class $\cald(\calc)$
      has the property (TREE) from Definition~\ref{def: class C_0}.

\item \label{the: Isomorphisms Conjecture and actions on trees: vcyc}
      Suppose  that $\calh^?_*$ has the regular tree property \eqref{def: tree-property} and
      $\calc \subseteq \calvcyc$. Then the class $\cald(\calc)$
      has the property (TREE$_\calr$) from Definition~\ref{def: class CL}.

\end{enumerate}
\end{theorem}
\begin{proof}
Let $G$ act on a tree $T$.
Denote by $V$ the set of vertices of $T$ and by $E$ the set of edges.
For $x \in V \cup E$ denote by $G_x$ the isotropy group of $x$ and by $\phi_x\colon G_x \to G$ the inclusion.
Let $\cali_T = \{ H \leq G \mid T^H \not= \emptyset \}$.
Since in a tree there is a unique geodesic between any two points,
the fixed set $T^H$ is contractible for $H \in \cali_T$.
Thus $T$ is a model for $\EGF{G}{\cali_T}$.

Next we prove  \ref{the: Isomorphisms Conjecture and actions on trees: fibered}.
In this case we assume that for each $x \in V \cup E$ the pair $(G_x,\calc(G_x))$ satisfies the
Fibered Isomorphism Conjecture~\ref{def: (Fibered) Isomorphism Conjectures for calh^?_*}.
Let $\phi\colon K \to G$ be a group homomorphism.
Then $K$ acts via $\phi$ on $T$.
Equipped with this action $T$ is also a model for $\EGF{K}{\phi^* \cali_T}$.
The tree property implies that $(K, \phi^* \cali_T)$ satisfies the
Isomorphism Conjecture~\ref{def: (Fibered) Isomorphism Conjectures for calh^?_*}.
Thus $(G,\cali_T)$ satisfies the
Fibered Isomorphism Conjecture~\ref{def: (Fibered) Isomorphism Conjectures for calh^?_*}.
Since the isotropy groups of $T$ satisfy the
Fibered Isomorphism Conjecture~\ref{def: (Fibered) Isomorphism Conjectures for calh^?_*}
with respect to $\calc$, we can deduce from Theorem~\ref{the: Reducing the family} that
$(G,\calc(G) \cap \cali_T)$ satisfies the
Fibered Isomorphism Conjecture~\ref{def: (Fibered) Isomorphism Conjectures for calh^?_*}.
Finally we use the fact that for the
\emph{Fibered} Isomorphism Conjecture~\ref{def: (Fibered) Isomorphism Conjectures for calh^?_*}
we can always enlarge the family (see \cite[Lemma 1.6]{Bartels-Lueck(2004)})
to conclude that the pair $(G,\calc(G))$ satisfies the
Fibered Isomorphism Conjecture~\ref{def: (Fibered) Isomorphism Conjectures for calh^?_*}.

Next we prove \ref{the: Isomorphisms Conjecture and actions on trees: fin}.
In this case we assume that for each $x \in V \cup E$ the pair $(G_x,\calc(G_x))$ satisfies the
Isomorphism Conjecture~\ref{def: (Fibered) Isomorphism Conjectures for calh^?_*}.
Arguing as above we conclude that $(G,\cali_T \cap \calc(G))$ satisfies the
Isomorphism Conjecture~\ref{def: (Fibered) Isomorphism Conjectures for calh^?_*}.
Finite groups cannot act without fixed points on trees
\cite[Theorem 15 in 6.1 on page 58 and 6.3.1 on page 60]{Serre(1980)}.
Therefore $\cali_T \cap \calc(G) = \calc(G)$.

Finally, we prove \ref{the: Isomorphisms Conjecture and actions on trees: vcyc}.
In this case we assume that for each $x \in V \cup E$ the pair $(G_x,\calc(G_x))$ satisfies the
Isomorphism Conjecture~\ref{def: (Fibered) Isomorphism Conjectures for calh^?_*}
and that $G_e$ is regular coherent for each $e \in E$.
Arguing as before we conclude that $(G,\cali_T \cap \calc(G))$ satisfies the
Isomorphism Conjecture~\ref{def: (Fibered) Isomorphism Conjectures for calh^?_*}.
We have to show that $(G,\calc(G))$ satisfies the
Isomorphism Conjecture~\ref{def: (Fibered) Isomorphism Conjectures for calh^?_*}.
Because of Theorem~\ref{the: Reducing the family} it suffices to show for any
virtually cyclic group $V \in \calc(G)$ that the
Isomorphism Conjecture~\ref{def: (Fibered) Isomorphism Conjectures for calh^?_*}
holds for $(V,\cali_T \cap \calc(G)|_V) = (V,\cali_T)$.

We first  consider the case, where $V$ contains a non-trivial normal finite subgroup $F$.
We saw above that $T^F$ is not empty and contractible.
By Lemma~\ref{lem: regular coherent groups are torsionfree} regular coherent groups are torsionfree.
Thus isotropy groups of edges are torsionfree, therefore $T^F$ is just a single vertex of $T$.
Since $F$ is normal in $V$, the action of $V$ leaves the fixed points of $F$ invariant.
Therefore the vertex $T^F$ is a fixed point for $V$. Hence we have $V \in \cali_T$
so that $\cali_T$ consists of all subgroups of $V$.

If $V$ does not contain a non-trivial normal finite subgroup $F$, it is either
$\IZ$ or the infinite dihedral group. In both cases $V$ acts on the tree $\IR$
with finite stabilizers such that the stabilizers of the edges are trivial
and every finite subgroup of $V$ occurs as stabilizer. The tree $\IR$ is a model
for $\EGF{V}{\calfin}$.
Since $\calh^?_*$ has the regular tree property \eqref{def: tree-property}
the map $\calh_n^V(\IR) \to \calh_n^V(\pt)$ is bijective for all $n \in \IZ$.
This shows that $V$ satisfies the Isomorphism Conjecture~\ref{def: (Fibered) Isomorphism Conjectures for calh^?_*}
for $(V,\calfin(V))$.
If $V = \IZ$, then  every subgroup $H \subseteq V$ is trivial
or isomorphic to $\IZ$. If $V$ is the infinite dihedral group, then
any subgroup $H$ of $V$ is finite, infinite cyclic or infinite dihedral.
We conclude from Theorem~\ref{the: Reducing the family}
that $V$ satisfies the Isomorphism Conjecture~\ref{def: (Fibered) Isomorphism Conjectures for calh^?_*}
for every family which contains $\calfin$, in particular for $\cali_T$.
\end{proof}

\begin{lemma} \label{lem: regular coherent groups are torsionfree}
Regular coherent groups are torsionfree.
\end{lemma}

\begin{proof}
Assume that $F$ is a finite subgroup of a regular coherent group $G$.
Then the $\IZ G$-module $\IZ[G/F]$ is finitely presented and has a finite-dimensional resolution
by finitely generated projective $\IZ G$-modules since $G$ is regular coherent and the ring $\IZ$ is regular.
Thus the restriction of $\IZ[G/F]$ to a $\IZ F$-module has a finite-dimensional resolution
by projective (but no longer finite generated) $\IZ F$-modules.
As an $\IZ F$-module $\IZ[G/F]$ contains $\IZ$ (with the constant $F$-action) as a direct summand.
Therefore $\IZ$ has a finite-dimensional resolution by projective $\IZ F$-modules.
This is only possible if $F$ is the trivial group.
\end{proof}

\begin{corollary}
\label{cor: extensions by groups acting on trees with finite stabilizers}
Let $\calh^?_*$ be an equivariant homology theory which has the tree
property (see Definition \ref{def: tree-property}).
Let $1 \to K \to G \to Q \to 1$ be an extension of groups.
Suppose that $K$ acts on a tree with finite stabilizers and that
$(Q,\calfin)$ satisfies the Fibered Isomorphism Conjecture
\ref{def: (Fibered) Isomorphism Conjectures for calh^?_*}
(in the range $\le N$).

Then $(G,\calfin)$ satisfies the Fibered Isomorphism Conjecture
\ref{def: (Fibered) Isomorphism Conjectures for calh^?_*}
(in the range $\le N$).
\end{corollary}
\begin{proof}
We first treat the case $Q = \{1\}$. Then the claim follows from
Theorem \ref{the: Isomorphisms Conjecture and actions on trees}
\ref{the: Isomorphisms Conjecture and actions on trees: fibered}
because for  a finite
group $F$ the pair $(F,\calfin)$ obviously satisfies the
Fibered Isomorphism Conjecture
\ref{def: (Fibered) Isomorphism Conjectures for calh^?_*}.

Next we treat the case, where $Q$ is finite. By a result of Dunwoody
\cite[Theorem 1.1]{Dunwoody(1979)} a group $K$
acts on a tree with finite stabilizers if and only
if $H^p(K;\IQ) = 0$ for each $p \ge 2$.
Since $K$ acts on a tree with finite
stabilizers, the trivial $\IQ K$-module $\IQ$ has a
$1$-dimensional projective resolution.
Hence  the trivial $\IQ G$-module $\IQ$ has a
$1$-dimensional projective resolution
since $[G:K]$ is finite and invertible in $\IQ$. This implies
$H^p(G;\IQ) = 0$ for each $p \ge 2$.
Hence also $G$ acts on a tree with finite stabilizers
if $Q$ is finite. This proves the claim  for finite $Q$.

Now the general case follows from
Lemma \ref{lem: fibered conjecture and extensions}.
\end{proof}

\begin{lemma}
\label{lem: suffices to check tree property on amal. prod. and HNN}
Let $\calh_*^?$ be an equivariant homology theory which is continuous.
Then the following assertions are equivalent.
\begin{enumerate}
\item
\label{lem: suffices to check tree property on amal. prod. and HNN: general}
For each $1$-dimensional $G$-$CW$-complex $T$
for which each component is contractible
(after forgetting the group action),
the projection $\pr_T\colon T \to \pi_0(T)$ induces isomorphisms
$$
\calh_n^G(\pr_T) \colon \calh_n^G(T) \xrightarrow{\cong}
           \calh_n^G(\pi_0(T)),
$$
for each $n \in \IZ$, where we consider $\pi_0(T)$ as a
$G$-space using the discrete topology;
\item
\label{lem: suffices to check tree property on amal. prod. and HNN: tree property}
$\calh^?_*$ has the tree property, i.e. for each $1$-dimensional
$G$-$CW$-complex $T$, which is contractible
(after forgetting the group action), and each $n \in \IZ$
we obtain isomorphisms
$$
\calh_n^G(\pr_T) \colon \calh_n^G(T) \xrightarrow{\cong}
               \calh_n^G(\pt);
$$
\item
\label{lem: suffices to check tree property on amal. prod. and HNN: few cells}
For each $1$-dimensional
$G$-$CW$-complex $X$, which is contractible
(after forgetting the group action) and has only one equivariant $1$-cell,
and each $n \in \IZ$ we obtain isomorphisms
$$
\calh_n^G(\pr_X) \colon \calh_n^G(X) \xrightarrow{\cong}
                \calh_n^G(\pt).
$$
\end{enumerate}
These three assertions remain equivalent if we add the requirement that
the isotropy groups of edges are regular coherent
to each assertion.
(Thus \ref{lem: suffices to check tree property on amal. prod. and HNN: tree property} becomes
the assertion that $\calh^?_*$ has the regular tree property.)
\end{lemma}
\begin{proof}
\ref{lem: suffices to check tree property on amal. prod. and HNN: general}
$\Rightarrow$
\ref{lem: suffices to check tree property on amal. prod. and HNN: tree property}
$\Rightarrow$
\ref{lem: suffices to check tree property on amal. prod. and HNN: few cells}
is obvious.
\\[2ex]
\ref{lem: suffices to check tree property on amal. prod. and HNN: few cells}
$\Rightarrow$
\ref{lem: suffices to check tree property on amal. prod. and HNN: general}
We prove the claim first under the assumption that $G\backslash T$
has finitely many $1$-cells.

We use induction over the number of $1$-cells in $G\backslash T$.
In the induction beginning, where $G\backslash T$ has no $1$-cell,
$T$ is the disjoint union of homogeneous spaces and the claim follows from the
fact that $\calh^?_*$ satisfies the disjoint union axiom.

In the induction step we can write $T$ as a $G$-pushout
$$\comsquare{G/H \times S^0}{q}{T_0}{}{}{G/H \times D^1}{}{T}$$
for a $G$-$CW$-subcomplex $T_0 \subseteq T$ such that
$G\backslash T_0$ has one $1$-cell less than $G\backslash T$.
Here $H$ is the isotropy group of the $1$-cell of $T$ that is 
not contained in $T_0$.
Since a connected subgraph of a tree is again a tree, each component
of $T_0$ is contractible. The induction hypothesis applies to $T_0$,
$G/H \times S^0$ and $G/H \times D^1$. Define $X$ to be the
$G$-pushout
$$
\comsquare{G/H \times S^0}{\pr_{T_0} \circ q}{\pi_0(T_0)}{}{}
          {G/H \times D^1}{}{X}
$$
The $G$-maps $\pr_{T_0} \colon T_0 \to \pi_0(T_0)$, $\id_{G/H \times S^0}$ and
$\id_{G/H \times D^1}$ are non-equivariant homotopy equivalences and induce
a $G$-map
$f \colon T \to X$ which is a non-equivariant homotopy equivalence since
$G/H \times S^0 \to G/H \times D^1$ is a cofibration. In particular
$X$ is a $1$-dimensional $G$-$CW$-complex whose components are contractible.
By a Mayer-Vietoris argument
and the Five-Lemma the map
$$\calh^G_n(f) \colon \calh^G_n(T) \xrightarrow{\cong} \calh^G_n(X)$$
is bijective for all $n \in \IZ$. The following diagram commutes
$$\comsquare{T}{\pr_T}{\pi_0(T)}{f}{\pi_0(f)}{X}{\pr_X}{\pi_0(X)}$$
Since the map $\pi_0(f)$ is bijective and hence a $G$-homeomorphism,
$\calh^G_n(\pi_0(f))$ is bijective for all $n \in \IZ$.
Recall that we have to show that $\calh^G_n(\pr_T)$ is bijective for
all $n \in \IZ$. Hence it suffices to show that
$\calh^G_n(\pr_X)$ is bijective for all $n \in \IZ$.
This follows from the fact that we can write $X$ as a
disjoint union of a $G$-$CW$-complex $Y$,
for which the
assumption \ref{lem: suffices to check tree property on amal. prod. and HNN: few cells}
applies, and a $0$-dimensional $G$-$CW$-complex $Z$, for which the
induction beginning applies,
and that $\calh^?_*$ satisfies the disjoint union axiom.

Next we treat the general case. Because $\calh^?_*$ satisfies the disjoint
union axiom, we can assume without loss
of generality that  $G\backslash T$ is connected. Since we can write
$T = G \times_H T'$ for a path component $T'$
and we have natural isomorphisms
$\calh^H(T') \xrightarrow{\cong} \calh_n^G(T)$ and
$\calh^H(\pt) \xrightarrow{\cong} \calh_n^G(G/H)$,
we can assume without loss of generality
that $T$ is contractible.

Fix a $0$-cell $e \in G\backslash T$.
Let $I$ be the set of finite connected $CW$-subcomplexes
$Z \subseteq G\backslash T$
with $e \in Z$.
It can be directed by inclusion and satisfies
$G\backslash T = \bigcup_{Z \in I} Z$.
Let $p \colon T \to G\backslash T$ be the projection.
Then $T$ is the directed union of the $G$-$CW$-subcomplexes $p^{-1}(Z)$.
Because of Lemma~\ref{lem: G-homology theory and colimit}
the canonical map
$$
\colim_{Z \in I} \calh^G_n(p^{-1}(Z)) \xrightarrow{\cong}
\calh^G_n(T)
$$
is bijective. Since each $G$-$CW$-complex $p^{-1}(Z)$ has only finitely
many equivariant $1$-cells and hence
satisfies the claim, and a colimit of a system of isomorphisms
is again an isomorphism, it suffices to show that
\begin{eqnarray}
\colim_{Z \in I} \calh^G_n(\pi_0(p^{-1}(Z)) & \to & \calh^G_n(\pt)
\label{colim over pi_0-s}
\end{eqnarray}
is bijective. Fix $\widetilde{e} \in T$ with $p(\widetilde{e}) = e$. Let
$G_Z$ be the isotropy group of the path component of $p^{-1}(Z)$ containing
$\widetilde{e}$ in the $G$-set
$\pi_0(p^{-1}(Z))$. Since each $Z$ is connected, $\pi_0(p^{-1}Z)$ is  $G/G_Z$.
We have for any inclusion $Z_1 \subseteq Z_2$ for
elements $Z_1, Z_2 \in I$, that $G_{Z_1}$ is a subgroup of
$G_{Z_2}$. We have $G = \bigcup_{Z \in I} G_Z$. Since $\calh^?_*$ is continuous,
we get an isomorphism
$$\colim_{Z \in I} \calh^G_n(G/G_Z) \xrightarrow{\cong}  \calh^G_n(\pt).$$
But this isomorphism can easily be identified with the map
\eqref{colim over pi_0-s}. This finishes the proof of
Lemma \ref{lem: suffices to check tree property on amal. prod. and HNN}.
\end{proof}

\begin{remark} \label{rem: properties (TREE) and (TREE')}
Let $G$ act on a tree $T$, such that $G \backslash T$ has only finitely many $1$-cells.
The proof of Lemma~\ref{lem: suffices to check tree property on amal. prod. and HNN}
shows that then $G$ acts on tree $X$ with the following properties:
The quotient $G \backslash X$ has only one $1$-cell.
For each edge $e$ of $X$ the isotropy group $G_e$ is also the isotropy group of an edge $e'$ of $T$.
For each vertex $v$ of $X$ there is a subtree $T_v$ of $T$ that is invariant under the isotropy group $G_v$
and for which $G_v \backslash T_v$ has one less $1$-cell than $G \backslash T$.
In combination with the colimit argument from the proof of
Lemma~\ref{lem: suffices to check tree property on amal. prod. and HNN} this means that a class of groups $\calc$
that has property (COL) from Definition~\ref{def: class C_0} has property (TREE) from Definition~\ref{def: class C_0}
if and only if it has the following property
\begin{description}
\item[(TREE')] Suppose that $G$ acts on a tree $T$ where $T$ has only one equivariant $1$-cell.
               Assume that for each $x \in T$ the isotropy group $G_x$ belongs to $\calc$.
               Then $G$ belongs to $\calc$;
\end{description}
and has property (TREE$_\calr$) if and only if it has the property
\begin{description}
\item[(TREE'$_\calr$)] Suppose that $G$ acts on a tree $T$.
               Assume that for each $x \in T$ the isotropy
               group $G_x$ belongs to $\calc$. For each edge $e$ of $T$, assume that
               the isotropy group $G_e$ is regular coherent.
               Then $G$ belongs to $\calc$.
\end{description}
Note on the other hand, that the statement that the
Fibered Isomorphism Conjecture~\ref{def: (Fibered) Isomorphism Conjectures for calh^?_*}
has property (TREE') is really a statement about arbitrary actions on trees:
If $G$ acts on a tree $T$ where $T$ has only one equivariant $1$-cell and $\phi \colon K \to G$
is a group homomorphism, then the induced action of $K$ on $T$ may have more equivariant $1$-cells and may
even be no longer cocompact.
Therefore we have to consider general trees in the formulation of the tree property
in Definition~\ref{def: tree-property}.
\end{remark}


\typeout{----- Equivariant homology theories constructed from spectra ------}
\section{Equivariant homology theories constructed from spectra}
\label{sec: Equivariant homology theories constructed from spectra}

In this section we want to
give a criterion when an equivariant homology theory has
the tree property provided that
it arises from  a covariant functor $\bfE \colon \GROUPOIDS \to \SPECTRA$
which sends equivalences of groupoids to weak equivalences of spectra.
This will be the main example for us.

Fix a group $G$. The \emph{transport groupoid} $\calg^G(S)$ of a $G$-set $S$
has $S$ as set of objects and the set of morphism from $s_1$ to $s_2$ consists
of those element $g \in G$ with $gs_1 = s_2$. Composition of
morphisms comes from the
group structure on $G$. The orbit category $\Or(G)$ has as objects
homogeneous spaces
$G/H$ and as morphisms $G$-maps. We obtain a covariant functor
$\calg^G \colon \Or(G) \to \GROUPOIDS, \hspace{1mm} G/H \mapsto \calg^G(G/H)$.
Define the covariant functor
$\bfE^G \colon \Or(G) \to \SPECTRA$ by $\bfE \circ \calg^G$.
Let $H^G_*(-;\bfE)$ be the
$G$-homology theory associated to $\bfE^G$ in
\cite[Section 4 and Section 7]{Davis-Lueck(1998)}. It is not hard to
construct the
relevant induction structure to get an equivariant homology theory
$H^?_*(-;\bfE)$.
It has the property that for each
group $G$ with subgroup $H \subseteq G$ and each $n \in \IZ$
we have canonical isomorphisms
$$H^G_n(G/H;\bfE) \cong H_n^H(\pt;\bfE) \cong \pi_n(\bfE(H)).$$
In the expression $\bfE(H)$ we think of the group $H$ as a groupoid
with one object.
More details of the construction of $H^?_*(-;\bfE)$ can be found in
\cite{Lueck-Reich(2004g)} and \cite{SauerJ(2002)}.

\begin{lemma} \label{lem: tree property for calh^?_*(-,bfE)}
The equivariant homology theory $H^?_*(-,\bfE)$ is continuous
and has the tree property if and only if the following  conditions are satisfied
\begin{enumerate}
\item \label{lem: tree property for calh^?_*(-,bfE): continuity}
For each group $G$  and directed system of subgroups $\{G_i \mid i \in I\}$,
which is directed by inclusion and satisfies $\bigcup_{i \in I} G_i = G$, and
each  $n \in \IZ$ the map
$$
\colim_{i \in I} j_i \colon \colim_{i \in I} \pi_n(\bfE(G_i))
                   \to \pi_n(\bfE(G))
$$
is an isomorphism, where $j_i$ is the homomorphism induced by the
inclusion $G_i \to G$;
\item  \label{lem: tree property for calh^?_*(-,bfE): amalgamated product}
Consider a pushout of groups
\begin{eqnarray}
\comsquare{H_0}{i_1}{H_1}{i_2}{j_1}{H_2}{j_2}{G}
\label{amalgamation pushout}
\end{eqnarray}
such that $i_1$ and $i_2$ are injective. In other words, $G$ is the
amalgamated product of $H_1$ and $H_2$ over $H_0$ with respect to the
injections $i_1$ and $i_2$. Then for each such
pushout \eqref{amalgamation pushout}
the following square of spectra is homotopy cocartesian
\begin{eqnarray}
\comsquarel{\bfE^G(G/H_0)\vee \bfE^G(G/H_0)}
          {\bfE^G(\pr_1)\vee \bfE^G(\pr_2)}
          {\bfE^G(G/H_1) \vee \bfE^G(G/H_2)}
          {\bfid \vee  \bfid}{\bfE^G(\pr_3) \vee \bfE^G(\pr_4)}
          {\bfE^G(G/H_0)}{\bfE^G(\pr_5)}{\bfE^G(G/G)}
\label{amalgamation pushout with bfE^G}
\end{eqnarray}
where the maps labeled $\pr_i$ denote canonical projections.

\item \label{lem: tree property for calh^?_*(-,bfE): HNN}
Let $i_0, i_1 \colon H \to K$ be injective group homomorphisms.
Let $G$ be the HNN-extension associated to $i_0$ and $i_1$.
The HNN-extension comes with an inclusion $j \colon K \to G$
and $t \in G$ such that $j \circ i_0 = c_t \circ j \circ i_1$,
where $c_t$ is conjugation by $t$. (This is the defining property
of the HNN-extension.)
We will use $i_0$ to consider $H$ as a subgroup of $G$.
Then the following
square of spectra is homotopy cocartesian
\begin{eqnarray}
\comsquare{\bfE^G(G/H) \vee \bfE^G(G/H)}
{\bfE^G(\pr_0) \vee \bfE^G(\beta)}
{\bfE^G(G/K)}
{\bfid \vee  \bfid}{\bfE^G(\pr_1)}
{\bfE^G(G/H)}{\bfE^G(\pr_2)}{\bfE^G(G/G)}
\label{HNN pushout with bfE^G}\end{eqnarray}
where the maps labeled $\pr_i$ are
canonical projections while
$\beta$ is defined by
$\beta(gH) = gtK$.
\end{enumerate}
For conditions \ref{lem: tree property for calh^?_*(-,bfE): amalgamated product}
and \ref{lem: tree property for calh^?_*(-,bfE): HNN}
one can also consider the regular versions 
where 
\eqref{amalgamation pushout with bfE^G}
is only required to be homotopy cartesian if in addition $H_0$ is
regular coherent and torsionfree
and 
\eqref{HNN pushout with bfE^G}
is only required to be homotopy cartesian if in addition $H$ is
regular coherent and torsionfree. 
Then the equivariant homology theory $H^?_*(-,\bfE)$ is continuous
and has the regular tree property if and only
the condition \ref{lem: tree property for calh^?_*(-,bfE): continuity}
and the regular versions of
conditions \ref{lem: tree property for calh^?_*(-,bfE): amalgamated product}
and \ref{lem: tree property for calh^?_*(-,bfE): HNN} hold.
\end{lemma}

\begin{proof}
Obviously condition \ref{lem: tree property for calh^?_*(-,bfE): continuity}
is equivalent to the condition
that $\calh^?_*$ is continuous. From now on we assume
that $\calh^?_*$ is continuous.

Suppose that the two conditions
\ref{lem: tree property for calh^?_*(-,bfE): amalgamated product} and
\ref{lem: tree property for calh^?_*(-,bfE): HNN} are satisfied. Because of
Lemma~\ref{lem: suffices to check tree property on amal. prod. and HNN}
it suffices to prove the tree
property only for $1$-dimensional contractible $G$-$CW$-complexes $T$ such that
there is precisely one equivariant $1$-cell.
Such a $G$-$CW$-complex will have precisely one or precisely two
equivariant $0$-cells.
We only treat the case, where there are two equivariant $0$-cells,
the proof of the other case is analogous using condition
\ref{lem: tree property for calh^?_*(-,bfE): HNN} instead of
condition \ref{lem: tree property for calh^?_*(-,bfE): amalgamated product}.

We can write $T$ as a $G$-pushout
$$
\comsquare{G/H_0 \times S^0}{\pr_1 \coprod \pr_2}
  {G/H_1 \coprod G/H_2}{}{}{G/H_0 \times D^1}{}{T}
$$
where $H_0$ is a subgroup of both $H_1$ and $H_2$ and $\pr_1$ and
$\pr_2$ are the canonical projections.
Recall that a $G$-space $Z$ defines a contravariant functor
$\Or(G) \to \SPACES, \; G/H \mapsto \map_G(G/H,Z)$
and that we get a spectrum
$\map_G(G/?,Z) \wedge_{\Or(G)} \bfE^G$ by the tensor product over
the orbit category (see \cite[Section 1]{Davis-Lueck(1998)}).
If we apply $\map_G(G/?,-) \wedge_{\Or(G)} \bfE^G$ to the $G$-pushout above,
we obtain a homotopy cocartesian diagram of spectra
\begin{eqnarray}
& \comsquare{\bfE^G(G/H_0) \vee \bfE^G(G/H_0)}
            {\bfE^G(\pr_1) \vee \bfE(\pr_2)}
            {\bfE^G(G/H_1) \vee \bfE^G(G/H_2)}
            {\bfid \vee  \bfid }{}
            {\bfE^G(G/H_0)}{}{\map_G(G/?,T)\wedge_{\Or(G)} \bfE^G}
& \label{hococadia I}
\end{eqnarray}
The following diagram is a pushout of groups
$$\comsquare{H_0}{i_1}{H_1}{i_2}{j_1}{H_2}{j_2}{G}$$
where $i_k \colon H_0 \to H_k$, $j_k \colon H_k \to G$ are inclusion
(see \cite[Example 1 on page 43]{Serre(1980)}).
Hence by condition
\ref{lem: tree property for calh^?_*(-,bfE): amalgamated product}
we have the homotopy cocartesian
square \eqref{amalgamation pushout with bfE^G}.
The projection $\pr \colon T \to G/G$
induces a map from the right lower corner of the diagram \eqref{hococadia I}
to the right lower corner of the diagram \eqref{amalgamation pushout with bfE^G}, if we identify
$G/G \wedge_{\Or(G)} \bfE^G = \bfE^G(G/G)$.
If we take the identity on the other three corners, we get a map between
homotopy cocartesian squares of spectra.
Since the three identity maps are obviously weak equivalences, the fourth map
induced by the projection is a weak equivalence.
But this map induces on homotopy groups
the map $H_n^G(\pr;\bfE) \colon H^G_n(T;\bfE) \to H^G_n(\pt)$ which
is hence bijective for each $n \in \IZ$.

This shows that $H^?_*(-;\bfE)$ has the tree property if conditions
\ref{lem: tree property for calh^?_*(-,bfE): amalgamated product} and
\ref{lem: tree property for calh^?_*(-,bfE): HNN} are satisfied.
It is now also obvious that
conditions \ref{lem: tree property for calh^?_*(-,bfE): amalgamated product}
and
\ref{lem: tree property for calh^?_*(-,bfE): HNN}
hold if $H^?_*(-;\bfE)$ has the tree property.
\end{proof}

\pagebreak[1]

\begin{remark} \
\label{rem: smaller diagrams}
\begin{enumerate}
\item
\label{rem: smaller diagrams:amalgamated}
In the situation of
Lemma~\ref{lem: tree property for calh^?_*(-,bfE)}~
\ref{lem: tree property for calh^?_*(-,bfE): amalgamated product}
diagram \eqref{amalgamation pushout with bfE^G}
is homotopy cocartesian if and only if
the commutative diagram
$$
\comsquare{\bfE(H)\vee \bfE(H)}{\bfE(i_1)\vee \bfE(i_2)}
          {\bfE(G_1) \vee \bfE(G_2)}
          {\bfid \vee  \bfid}{\bfE(j_1) \vee \bfE(j_2)}
          {\bfE(H)}{\bfE(j_0)}{\bfE(G)}
$$
where $j_0 \colon H \to G$ is defined to be
$j_1 \circ i_1 = j_2 \circ i_2$,
is homotopy cocartesian since there is a canonical
weak equivalences from each  corner of
this square to the corresponding corner of
\eqref{amalgamation pushout with bfE^G}.
\item
\label{rem: smaller diagrams:HNN}
The situation in
Lemma~\ref{lem: tree property for calh^?_*(-,bfE)}~
\ref{lem: tree property for calh^?_*(-,bfE): HNN}
is a bit more complicated.
The natural diagram to consider is
\begin{eqnarray}
\comsquare{\bfE(H) \vee \bfE(H)}
  {\bfE(i_0) \vee \bfE(i_1)}{\bfE(K)}{\bfid \vee  \bfid}{\bfE(j)}
  {\bfE(H)}{\bfE(j \circ i_0)}{\bfE(G)}
\label{HNN pushout}
\end{eqnarray}
However,
\eqref{HNN pushout} is not commutative, while
\eqref{HNN pushout with bfE^G} is commutative.
There is a canonical weak equivalence from each corner of
\eqref{HNN pushout with bfE^G}
to the corresponding corner of \eqref{HNN pushout},
but those maps do not make the square
\[
\comsquare{\bfE(H) \vee \bfE(H)}{\bfE(i_0) \vee \bfE(i_1)}{\bfE(K)}
          {}{}
{\bfE^G(G/H) \vee \bfE^G(G/H)}
{E^G(\pr_0) \vee E^G(\beta)}
{\bfE^G(G/K)}
\]
commutative.

The failure of the commutativity of \eqref{HNN pushout}
stems from the fact, that the underlying diagram of groups
commutes only up to conjugation, i.e.
$j \circ i_0 \neq j \circ i_1 = c_t \circ j \circ i_0$.
It is a consequence of the definitions
that $\bfE(c_t)$ is weakly homotopic to $\id_{\bfE(G)}$, but in general
there is no preferred homotopy. On the other hand
$\bfE\colon \GROUPOIDS \to \SPECTRA$ is often slightly better than required
in the discussion before
Lemma~\ref{lem: tree property for calh^?_*(-,bfE)},
namely $\bfE$ is a $2$-functor. This means that
if $\tau$ is a natural transformation
between functors $f,g$ between groupoids, then there is a (preferred) homotopy
$\bfE(\tau)$ from $\bfE(f)$ to $\bfE(g)$.
Under this stronger assumption on $\bfE$
there is a canonical homotopy that makes \eqref{HNN pushout}
homotopy commutative and then condition
\ref{lem: tree property for calh^?_*(-,bfE): HNN}
in Lemma~\ref{lem: tree property for calh^?_*(-,bfE)}
is equivalent to requiring
that \eqref{HNN pushout} is homotopy cocartesian with respect to
the canonical  homotopy.
\end{enumerate}
\end{remark}


\typeout{--------- Isomorphisms Conjecture for spectra ----------------}
\section{Isomorphism Conjectures for spectra}
\label{sec: Isomorphism Conjectures for spectra}

In this section we relate the (Fibered) Isomorphism Conjecture
\ref{def: (Fibered) Isomorphism Conjectures for calh^?_*}
for an  equivariant homology theory $\calh^?_*$ to the versions appearing
in Farrell-Jones \cite{Farrell-Jones(1993a)} for algebraic
$K$- and $L$-theory.

Consider a group homomorphism $\phi \colon K \to G$,
a $K$-$CW$-complex $Z$ and a
covariant functor $\bfE \colon \SPACES \to \SPECTRA$,
which sends weak equivalences to
weak equivalences and is \emph{compatible with disjoint unions}, i.e.
for a family $\{Y_i \mid i \in I\}$ of spaces
the map induced by the inclusions $j_i \colon Y_i \to \coprod_{i \in I} Y_i$
$$
\bigvee_{i \in I} \bfE(j_i) \colon \bigvee_{i \in I} \bfE(Y_i)
\to \bfE\left(\coprod_{i \in I} Y_i\right)
$$
is a weak equivalence. We obtain a covariant functor
$$
\bfE^K_Z \colon \Or(K) \to \SPECTRA, \quad K/H \mapsto \bfE(Z \times_K K/H).
$$
Recall that for each covariant functor $\bfF \colon \Or(K) \to \SPECTRA$
there is a $K$-homology theory $H^K_*(-;\bfF)$
defined for $K$-$CW$-complexes with the property that $H_n^K(K/H;\bfF) \cong
\pi_n(\bfF(K/H))$ holds for $H \subseteq K$ and $n \in \IZ$.
\cite[Section 4 and Section 7]{Davis-Lueck(1998)}.
We denote by $\phi_*Z$ the $G$-space $G \times_{\phi} Z$ obtained by
induction with $\phi$ from the $K$-space $Z$. For a $G$-space $X$
let $\phi^*X$ be the
$K$-space obtained by restricting the $G$-action to a $K$-action using $\phi$.

\begin{lemma} \label{lem: descend}
For any $G$-$CW$-complex $X$ there is an isomorphism,
natural in $X$, $Z$ and $\bfE$,
$$
\phi_* \colon H_n^K(\phi^* X;\bfE^K_Z) \xrightarrow{\cong}
           H_n^G(X;\bfE^G_{\phi_*Z}).
$$
\end{lemma}
\begin{proof}
Let $\phi \colon \Or(K) \to \Or(G), \quad K/H \to G/\phi(H)$ be the
functor induced by
$\phi$. Given a contravariant (pointed) $\Or(G)$-space $A$ and
a covariant (pointed) $\Or(K)$-space $B$ there is an adjunction
\begin{eqnarray}
\res_{\phi} A \otimes_{\Or(K)} B & \xrightarrow{\cong} &
                              A \otimes_{\Or(G)} \ind_{\phi} B,
\label{adjunction ind and res}
\end{eqnarray}
where $\res_{\phi}$ is restriction and $\ind_{\phi}$
denotes induction with the functor
$\phi \colon \Or(K) \to \Or(G)$ (see \cite[Lemma 1.9]{Davis-Lueck(1998)}).
It induces a natural isomorphism
$$H_n^K(\phi^* X;\bfE^K_Z) \xrightarrow{\cong} H_n^G(X;\ind_{\phi}\bfE^K_Z)$$
There is a weak equivalence of covariant $\Or(G)$-spectra
$$\ind_{\phi}\bfE^K_Z \xrightarrow{ \cong} \bfE^G_{\phi_* Z}$$
coming from
\begin{multline*}
$$
\map_G(G/\phi(?),G/??))  \otimes_{\Or(K)} (Z \times_K K/?)
      \xrightarrow{\cong}  Z \times_{\phi} G/??,
\\
\quad (f,(z,k?)) \mapsto (z,f(\phi(k)\phi(?)))
\end{multline*}
and the fact that $\bfE$ is compatible with disjoint unions.
\end{proof}

\begin{lemma} \label{lem: reformulation of the fibered conjecture}
Let $\calf$ be a family of subgroups of $G$. Let $N \in \IZ$.
Then the following assertions are equivalent:

\begin{enumerate}

\item
\label{lem: reformulation of the fibered conjecture: Non fibered, non-connected}
For any free  $G$-$CW$-complex $Z$ and $n \in \IZ$ (with $n \leq N$) the assembly map
$$H_n^G(\EGF{G}{\calf};\bfE^G_Z) \xrightarrow{\cong} H_n^G(\pt;\bfE^G_Z)$$
is bijective;

\item \label{lem: reformulation of the fibered conjecture: connected}
For each injective group homomorphism $\phi \colon K \to G$  and any
free connected $K$-$CW$-complex $Z$ and $n \in \IZ$ (with $n \leq N$) the assembly map
$$
H_n^K(\EGF{K}{\phi^*(\calf)};\bfE^K_Z)
  \xrightarrow{\cong} H_n^K(\pt;\bfE^K_Z)
$$
is bijective;

\item \label{lem: reformulation of the fibered conjecture: simply connected}
For each group homomorphism $\phi \colon K \to G$  and any free simply
connected $K$-$CW$-complex $Z$ and $n \in \IZ$ (with $n \leq N$) the assembly map
$$
H_n^K(\EGF{K}{\phi^*(\calf)};\bfE^K_Z) \xrightarrow{\cong}
        H_n^K(\pt;\bfE^K_Z)
$$
is bijective.

\end{enumerate}
\end{lemma}
\begin{proof}
\ref{lem: reformulation of the fibered conjecture: Non fibered, non-connected}
  $\Rightarrow$
\ref{lem: reformulation of the fibered conjecture: connected}
and
\ref{lem: reformulation of the fibered conjecture: Non fibered, non-connected}
 $\Rightarrow$
\ref{lem: reformulation of the fibered conjecture: simply connected}
These implications follow from   Lemma \ref{lem: descend}  since
for any group homomorphism $\phi \colon K \to G$
we have $\phi^*(\EGF{G}{\calf}) = \EGF{K}{\phi^*(\calf)}$.
\\[1mm]
\ref{lem: reformulation of the fibered conjecture: connected}  $\Rightarrow$
\ref{lem: reformulation of the fibered conjecture: Non fibered, non-connected}
We can write a $G$-$CW$-complex $Z$ as
 $Z = \coprod_{i \in I} G \times_{G_i} Z_i$ for subgroups $G_i \subseteq G$
and connected free $G_i$-$CW$-complexes $Z_i$. Since $\bfE$ is compatible
with disjoint unions, we conclude from \cite[Lemma 4.6]{Davis-Lueck(1998)}
that we can assume without loss of generality that $I$ consists of one
element $0$, i.e. $Z = G \times_{G_0} Z_0$. Now the claim follows from
Lemma \ref{lem: descend} applied to the inclusion $\phi \colon G_0 \to G$ and
the free connected $G_0$-$CW$-complex $Z_0$.
\\[1mm]
\ref{lem: reformulation of the fibered conjecture: simply connected}
$\Rightarrow$
\ref{lem: reformulation of the fibered conjecture: connected}
There is an extension of groups
$1 \to \pi_1(Z) \to \widetilde{K} \xrightarrow{p} K \to 1$
and a $\widetilde{K}$-action on the universal covering $\widetilde{Z}$
which extends
the $\pi_1(Z)$-action on $\widetilde{Z}$ and covers the $K$-action on $Z$.
Moreover, $\widetilde{Z}$ inherits the structure of a free
$\widetilde{K}$-$CW$-complex.
Now the claim follows from  Lemma \ref{lem: descend} applied to the group
homomorphism $p \colon \widetilde{K} \to K$ and the simply connected free
$\widetilde{K}$-$CW$-complex
$\widetilde{Z}$ since $p_*\widetilde{Z} = Z$ and
$p^*(\phi^* \calf) = (\phi \circ p)^*\calf$.
\end{proof}

\begin{lemma} \label{lem: bfE sends 2-connected maps to weak equivalences}
Suppose that for any two-connected map $f \colon X \to Y$ the induced map
$\bfE(f) \colon \bfE(X) \to \bfE(Y)$ is a weak equivalence. Let
$Z$ be a simply-connected $G$-$CW$-complex and $f \colon Z \to EG$ be the
classifying map.

Then it induces a weak equivalence of $\Or(G)$-spectra
$$\bff \colon \bfE^G_Z \to \bfE^G_{EG}$$
and in particular for each  $G$-$CW$-complex and each $n \in \IZ$ a natural
isomorphism
$$H_n^G(X;\bfE^G_Z) \xrightarrow{\cong} H_n^G(X;\bfE^G_{EG}).$$
\end{lemma}
\begin{proof}
The map $f \colon Z \to EG$ is $2$-connected.
Hence the induced map
$f \times_G \id_{G/H} \colon Z \times_G G/H \to EG \times_G G/H$ is
$2$-connected for all subgroups $H \subseteq G$.  Now apply
\cite[Lemma 4.6]{Davis-Lueck(1998)}.
\end{proof}

\begin{definition}[Fibered Isomorphism Conjecture for spectra]
\label{def: (Fibered) Isomorphisms Conjecture for spectra}
We say that a group $G$ satisfies the \emph{Isomorphism Conjecture for $\calf$
and $\bfE$ (in the range $\le N$)} if the assembly map
induced by the projection $\pr \colon \EGF{G}{\calf} \to \pt$
$$\asmb \colon H^G_n(\EGF{G}{\calf};\bfE^G_{EG})) \to H^G_n(\pt;\bfE^G_{EG}))$$
is bijective for all $n \in \IZ$ (with $n \le N$).

We say that a group $G$ satisfies the
\emph{Fibered Isomorphism Conjecture for $\calf$
and $\bfE$ (in the range $\le N$)}
if for any free $G$-$CW$-complex $Z$ the assembly map
induced by the projection $\pr \colon \EGF{G}{\calf} \to \pt$
$$\asmb \colon H^G_n(\EGF{G}{\calf};\bfE^G_Z)) \to
H^G_n(\pt;\bfE^G_Z))$$
is bijective  for all $n \in \IZ$ (with $n \le N$).
\end{definition}

\begin{remark}
\label{rem: relation to the fibered isomorphism conjecture of Farrell-Jones}
The (Fibered) Isomorphism Conjecture of
Farrell-Jones \cite{Farrell-Jones(1993a)}  for algebraic
$K$-theory or $L$-theory respectively is equivalent to
the (Fibered) Isomorphisms Conjecture
\ref{def: (Fibered) Isomorphisms Conjecture for spectra}
for $(G,\calvcyc,\bfE)$ if $\bfE \colon \SPACES \to \SPECTRA$ sends $X$ to the
non-connective algebraic $K$-theory spectrum or
$L^{\langle - \infty \rangle}$-theory spectrum
of the fundamental groupoid of $X$ respectively. Since a $2$-connected map
$f \colon X \to Y$
induces an equivalence on the fundamental groupoids,
Lemma  \ref{lem: bfE sends 2-connected maps to weak equivalences} applies.
Let $R$ be a ring. Consider the covariant functors
\begin{eqnarray*}
\bfK_R \colon \GROUPOIDS & \to & \SPECTRA\\
\bfL_R \colon \GROUPOIDS & \to & \SPECTRA
\end{eqnarray*}
defined in \cite[Section 2]{Davis-Lueck(1998)}  satisfying
$\pi_n(\bfK_R(G)) = K_n(RG)$ and
$\pi_n(\bfL_R(G)) = L_n^{\langle - \infty \rangle}(RG)$ for
each group $G$ and $n \in \IZ$.
Let $H^?_*(-,\bfK_R)$ and $H^?_*(-,\bfL_R)$ be the associated equivariant
homology theories. Then the (Fibered) Isomorphism Conjecture for
algebraic $K$-theory or algebraic $L$-theory
respectively for the group $G$ in the sense of
Farrell-Jones \cite{Farrell-Jones(1993a)}
is equivalent to the (Fibered) Isomorphisms Conjecture
\ref{def: (Fibered) Isomorphism Conjectures for calh^?_*} for
$H^?_*(-,\bfK_R)$ and $H^?_*(-,\bfL_R)$ for the pair $(G,\calvcyc)$.
This follows from Lemma
\ref{lem: reformulation of the fibered conjecture} and
Lemma  \ref{lem: bfE sends 2-connected maps to weak equivalences} .

For more information about the various conjectures such as the version for
pseudoisotopy or the Baum-Connes Conjecture we refer for instance to
\cite{Lueck-Reich(2004g)}.
\end{remark}


\typeout{---------- The $\KH$-Isomorphism Conjecture -------------}
\section{The KH-Isomorphism Conjecture}
\label{sec: The KH-Isomorphism Conjecture}

In this section we will formulate the $\KH$-Isomorphism Conjecture.
The construction of homotopy algebraic $K$-Theory
is a simplicial construction, so we will quickly fix the notation.
The category $\Delta$ has as objects finite ordered sets of the form
$\ULn = \{ 0 < 1 < \dots < n \}$ and order preserving maps as morphisms.
The n-simplex $\Delta^n_\bullet$ is the simplicial set
$\ULm \mapsto \Map_\Delta(\ULm,\ULn)$.
Let $R$ be a ring. The simplicial ring $R[\bullet]$ is defined by
\[
R[\ULn] = R[t_0,\dots,t_n] / (t_0 + \dots + t_n = 1).
\]
Here the structure maps acts as follows: if $f\colon \ULn \to \ULm$ is
order preserving then $f^*\colon R[\ULm] \to R[\ULn]$ is defined by
\[
f^*(t_k) = \sum_{j \in f^{-1}(k)} t_j.
\]
In \cite{Weibel(1989)} the homotopy algebraic $K$-theory $\KH_*(R)$ of $R$ is
defined as the homotopy groups of the realization $\bfKH R$ of the simplicial
spectrum $\bfK^{-\infty} R[\bullet]$. Here $\bfK^{-\infty}$ denotes the non-connected
$K$-theory spectrum; a construction is reviewed before
Definition~\ref{def: non-connective-Nil}.
To illustrate the construction of homotopy algebraic $K$-theory we give a proof of the
following fundamental property of homotopy algebraic $K$-theory, cf.~\cite[1.2.(i)]{Weibel(1989)}.

\begin{proposition} \label{pro: homotopy invariance}
The inclusion $R \hra R[X]$ gives an isomorphism
$\KH_n(R) \cong \KH_n(R[X])$ for all $n \in \IZ$.
\end{proposition}

\begin{proof}
It suffices to show that $R[\bullet] \to R[X][\bullet]$ is a homotopy
equivalence of simplicial rings (cf.\ Remark~\ref{ringhomotopy}).
To see this we need to show that
$X \mapsto 0$ is homotopic to the identity of $R[X][\bullet]$. Such a homotopy
$R[X][\bullet] \x \Delta^1_\bullet \to R[X][\bullet]$ is given by
\[
(X,f) \mapsto \left(\sum_{j \in f^{-1}(0)} t_j\right) X,
\]
where $f\colon \ULn \mapsto \UL{1}$.
\end{proof}

\begin{remark} \label{ringhomotopy}
If $S$ is a ring and $F$ is a finite set, then $S \x F$ has a ring structure
($S \x F \cong \prod_{f \in F} S$). Therefore we may view
$R[X][\bullet] \x \Delta^1_\bullet$ as a simplicial ring and the above homotopy
as a map of simplicial rings. Therefore we get a map
$| \bfK^{-\infty}(R[X][\bullet] \x \Delta^1_\bullet)| \to |\bfK^{-\infty} (R[X][\bullet])|$.
On the other hand, there is a map of simplicial spectra
$\bfK^{-\infty}(R[X][\bullet]) \x \Delta^1_\bullet \to
\bfK^{-\infty}(R[X][\bullet] \x \Delta^1_\bullet)$ defined as follows.
For $f \in \Delta^1_\ULn$ there is an obvious map of rings
$\iota_f\colon R[X][\ULn] \to R[X][\ULn] \x \{ f \} \hra R[X][\ULn] \x \Delta^1_n$.
Thus we can map $(x,f) \in \bfK^{-\infty}(R[X][\ULn]) \x \Delta^1_n$ to
$\bfK^{-\infty}(\iota_f)(x) \in \bfK^{-\infty}(R[X][\ULn] \x \Delta^1_n)$.
\end{remark}

In order to define an equivariant homology theory we
define the functor
\[
\bfKH_R\colon \GROUPOIDS \to \SPECTRA
\]
as the realization of the  simplicial functor
\[
\bfK_{R[\bullet]}\colon \GROUPOIDS \to \SPECTRA.
\]
Thus $\pi_n(\bfKH_R(G)) = \KH_n(RG)$.
Since the realization of a weak equivalence is again a weak equivalence,
$\bfKH_R$ sends equivalences of groupoids to weak equivalences of
spectra.

\begin{conjecture}[(Fibered) $\KH$-Isomorphism Conjecture]
\label{con: KH-Isomorphism-Conjecture}
A group $G$ is said to satisfy the (Fibered) $\KH$-Isomorphism Conjecture
(for a ring $R$)
if the pair $(G,\calfin)$ satisfies the
(Fibered) Isomorphism
Conjecture~\ref{def: (Fibered) Isomorphism Conjectures for calh^?_*}
for the equivariant homology theory $H^?_*(-;\bfKH_R)$.
\end{conjecture}

\begin{remark} \label{rem: VCYC or FIN in KH-Conjecture}
All virtually cyclic groups act on trees with finite stabilizers.
For a finite group $F$ the (Fibered) $\KH$-Isomorphism Conjecture
holds  (since $\EGF{F}{\calfin} = F/F$). Thus by
Theorem~\ref{the: inhertitance for KH} the (Fibered) $\KH$-Isomorphism Conjecture
holds for virtually cyclic groups.
Therefore Theorem~\ref{the: Reducing the family}
implies, that it makes no difference if we replace the family
of finite groups with the family of virtually cyclic groups in
the formulation of the (Fibered) $\KH$-Isomorphism Conjecture.
\end{remark}


\typeout{---------- The Relation between the K- and the KH-Isomorphism Conjecture -------------}
\section{The Relation between the $K$- and the $KH$-Isomorphism Conjecture}
\label{sec: The Relation between the K- and the KH-Isomorphism Conjecture}

There is a natural map $\bfK^{-\infty} R \to \bfKH R$ induced from the
inclusion of the constant simplicial ring $R$ into $R[\bullet]$.
Similarly we obtain a natural transformation $\bfK_R \to \bfKH_R$ of functors
from $\GROUPOIDS$ to $\SPECTRA$.
Thus we obtain a natural transformation of
equivariant homology theories
$H^?_*(-;\bfK_R) \to H^?_*(-;\bfKH_R)$ and a commutative diagram
between assembly maps
\begin{eqnarray}
& \comsquare{H^G_n(\EGF{G}{\calfin};\bfK_R)}{}{K_n(RG)}{}{}
{H^G_n(\EGF{G}{\calfin};\bfKH_R)}{}{\KH_n(RG)}&
\label{comparing-K-assembly-to-KH-assembly}
\end{eqnarray}
We first explain what the $\KH$-Isomorphism Conjecture~\ref{con: KH-Isomorphism-Conjecture}
implies for the $K$-Isomorphism Conjecture, i.e. the Farrell-Jones Conjecture for algebraic $K$-theory
(see Remark~\ref{rem: relation to the fibered isomorphism conjecture of Farrell-Jones}).
In order to state the  connection we need to recall the groups $N^p K_n(R)$
\cite[XIII]{Bass(1968)}.
They can be defined by $N^0 K_n(R)=K_n(R)$ and
\[
N^{p}K_n(R) = \Ker(N^{p-1}(q) \colon  N^{p-1}K_n(R[t]) \to N^{p-1}K_n(R)),
\]
where $q(t) = 0$.
For regular rings $N^{p}K_n(R) = 0$ for $p \geq 1$, see \cite{Gersten(1973)}.

\begin{proposition}
\label{pro: comparing-K-assembly-to-KH-assembly}
Let $G$ be a group that satisfies the
$\KH$-Isomorphism Conjecture~\ref{con: KH-Isomorphism-Conjecture} for the ring $R$.
\begin{enumerate}
\item \label{pro: comparing-K-assembly-to-KH-assembly:zero}
      Suppose that $N^{p}K_n(RF) = 0$ for all
      finite subgroups $F$ of $G$ and all $n \in \IZ$, $p \geq 1$.
      Then the assembly map with respect to the family $\calfin$
      in algebraic $K$-theory, i.e.\ the top row
      in \eqref{comparing-K-assembly-to-KH-assembly},
      is split injective.
\item \label{pro: comparing-K-assembly-to-KH-assembly:torsion}
      Suppose that $N^{p}K_n(RF) \ox \IQ = 0$ for all
      finite subgroups $F$ of $G$ and all $n \in \IZ$, $p \geq 1$.
      Then assembly map with respect to the family $\calfin$
      in algebraic $K$-theory, i.e.\ the top row
      in \eqref{comparing-K-assembly-to-KH-assembly}, is rationally split injective.
\end{enumerate}
\end{proposition}
\begin{proof}
By the spectral sequence from \cite[1.3]{Weibel(1989)} the canonical
map $K_*(A) \to \KH_*(A)$ is an isomorphism if $N^{p}K_n(A) = 0$
for all $n \in \IZ$ and $p \geq 1$ and a rational isomorphism
if $N^{p}K_n(A) \ox \IQ = 0$ for all $n \in \IZ$ and $p \geq 1$.
Therefore these assumptions imply by a spectral sequence argument
that the left vertical map in \eqref{comparing-K-assembly-to-KH-assembly} is an isomorphism
or a rational isomorphism respectively.
\end{proof}

\begin{remark} \label{rem: assumption of propositions hold in cases}
The assumptions of
Proposition~\ref{pro: comparing-K-assembly-to-KH-assembly}
\ref{pro: comparing-K-assembly-to-KH-assembly:zero} and
\ref{pro: comparing-K-assembly-to-KH-assembly:torsion}
are satisfied in the following cases.
\begin{enumerate}
\item If $R$ is a regular ring containing $\IQ$
      then $RF$ is regular
      for all finite groups $F$.
      Thus the assumption in
      \ref{pro: comparing-K-assembly-to-KH-assembly}
      \ref{pro: comparing-K-assembly-to-KH-assembly:zero}
      is satisfied.
\item If $R=\IZ$ then the assumption in
      \ref{pro: comparing-K-assembly-to-KH-assembly}
      \ref{pro: comparing-K-assembly-to-KH-assembly:torsion}
      is satisfied. This follows from \cite[6.4]{Weibel(1981)},
      which implies
      \[
      \NK_*(\IZ[t_1,\dots,t_n]F) \ox \IQ \cong
      \NK_*(\IQ[t_1,\dots,t_n]F).
      \]
      Thus for a finite group $F$ it follows that,
      $\NK_*(\IZ[t_1,\dots,t_n]F)$
      vanishes rationally, since $\IQ[t_1,\dots,t_n]F$ is regular.
      A straight forward induction shows that this implies
      $N^{p}K_n(\IZ F) \ox \IQ = 0$ for $p \geq 1$.
\item If $G$ is torsionfree and $R$ is regular, then the assumption in
      \ref{pro: comparing-K-assembly-to-KH-assembly}
      \ref{pro: comparing-K-assembly-to-KH-assembly:zero}
      is satisfied. In this case the Farrell-Jones Conjecture
      in algebraic $K$-theory asserts, that the top vertical map in
      \eqref{comparing-K-assembly-to-KH-assembly} is an isomorphism.
      Thus, in this situation ($R$ regular, $G$ torsionfree)
      the Farrell-Jones Conjecture in algebraic $K$-theory holds if
      $N^p_n(RG) = 0$ for all $n \in \IZ$ and
      $p \geq 1$ and if $G$ satisfies the $\KH$-Isomorphism Conjecture.
      \end{enumerate}
\end{remark}

Next we explain the reverse connection.

\begin{theorem}[The $K$-theory version implies the $\KH$-version]
\label{the: K-theory version implies KH-version}
Let $G$ be a group and let $R$ be a ring.
\begin{enumerate}
\item \label{the: K-theory version implies KH-version (polynom)}
      Suppose that the
      (Fibered) Farrell-Jones Conjecture in algebraic $K$-theory is true for
      $(G,R[x_1,x_2, \ldots x_n])$ for all $n \ge 1$
      then the (Fibered) $\KH$-Isomorphism
      Conjecture~\ref{con: KH-Isomorphism-Conjecture} is true for $(G,R)$.
\item \label{the: K-theory version implies KH-version (Laurent)}
      Suppose that the Fibered Farrell-Jones Conjecture in algebraic $K$-theory is true for
      $(G \times \IZ^n,R)$ for all $n \ge 1$. Then the Fibered $\KH$-Isomorphism
      Conjecture~\ref{con: KH-Isomorphism-Conjecture} is true for $(G,R)$.
\end{enumerate}
\end{theorem}

\begin{proof}
Let $\phi\colon K \to G$ be a group homomorphism.
The assembly map $H_n^K(\EGF{K}{\phi^*\calvcyc};\bfK_{R[\ULn]}) \to K_n(R[\ULn][K])$
is on the level of spectra given by the map
\[
\EGF{K}{\phi^*\calvcyc} \ox_{\Or(K)} \bfK_{R[\ULn]} \to K/K \ox_{\Or K} \bfK_{R[\ULn]} \simeq \bfK_{R[\ULn]}(K)
\]
induced by $\EGF{K}{\phi^*\calvcyc} \to K/K$.
The assumption in \ref{the: K-theory version implies KH-version (polynom)}
is that this map of spectra is a weak equivalence.
Using the fact that the realization of a map of simplicial spectra
that is levelwise a weak equivalences is a weak equivalence and the identification
\[
|\EGF{K}{\phi^*\calvcyc} \ox_{\Or(K)} \bfK_{R[\bullet]}| \cong \EGF{K}{\calvcyc} \ox_{\Or(K)} |\bfK_{R[\bullet]}|
\]
we conclude that the (Fibered) Farrell-Jones Conjecture for $(G, R[x_1,\dots,x_n])$ for all $n$
implies the (Fibered) $\KH$-Isomorphism Conjecture for $(G, R)$
with the family of finite subgroups replaced by the family of virtually cyclic subgroups.
By Remark~\ref{rem: VCYC or FIN in KH-Conjecture} this is equivalent to
the (Fibered) $\KH$-Isomorphism Conjecture~\ref{con: KH-Isomorphism-Conjecture} for $(G, R)$.

Next we prove \ref{the: K-theory version implies KH-version (Laurent)} by reducing it to
\ref{the: K-theory version implies KH-version (polynom)}.
For a group $K$ we denote by $p_K\colon K \x \IZ \to K$ the canonical projection.
We observe first that the (Fibered) Isomorphism Conjecture for $(G,\calvcyc,R[\IZ])$ is equivalent
to the (Fibered) Isomorphism Conjecture $(G \x \IZ, (p_G)^*\calvcyc, R)$ because
for every group $K$ and every $K$-space $X$
there is a natural isomorphism
\[
H_n^K( X; \bfK_{R[\IZ]}) \cong H_n^{K \x \IZ}( {p_K}^* X; \bfK_{R} )
\]
where ${p_K}^* X$ denotes the $K \x \IZ$-space obtained by restriction of $X$ along $p_K$
and because ${p_K}^* \EGF{K}{\calf} = \EGF{K \x \IZ}{(p_K)^*\calf}$ for every family of subgroups of $K$.
If the {\em Fibered} Isomorphism Conjecture holds for a family $\calf$, then it will also hold
for every family $\calg$ that contains $\calf$ \cite[Lemma 1.6]{Bartels-Lueck(2004)}.
Because the family of virtually cyclic subgroups of $G \x \IZ$ is contained in $(p_G)^*\calvcyc$
the Fibered Farrell-Jones Conjecture for $(G \x \IZ, R)$ implies the
Fibered Farrell-Jones Conjecture for $(G, R[\IZ])$.
By the Bass-Heller-Swan splittings \cite{Grayson(1976)}
the latter is equivalent to the Fibered Farrell-Jones Conjecture for $(G, R[x])$.
By induction on $n$ this means that the assumption in \ref{the: K-theory version implies KH-version (Laurent)}
implies the assumption of \ref{the: K-theory version implies KH-version (polynom)}.
\end{proof}

\begin{remark}
It is not unreasonable to expect that the non-fibered version of
Theorem~\ref{the: K-theory version implies KH-version}~\ref{the: K-theory version implies KH-version (Laurent)}
is also valid.
Our argument  would also prove the non-fibered version if we were to know that for every virtually cyclic
group $V$ the product $V \x \IZ$ satisfies the Farrell-Jones Conjecture.
This seems very likely, but we could not find such a statement in the literature.
\end{remark}

\begin{remark} \label{rem: consequences of the conjectures}
Let us briefly list some consequences of the Farrell-Jones Conjecture for algebraic
$K$-theory. Suppose that the Farrell-Jones Conjecture for algebraic
$K$-theory holds for the group $G$ and every regular ring $R$. Now consider a group $G$ and a regular ring $R$
with the property that either $\IQ \subseteq R$ holds or  $G$ is torsionfree.
The proof of
Proposition~\ref{pro:  comparing-K-assembly-to-KH-assembly},
Remark~\ref{rem: assumption of propositions hold in cases} and
Theorem~\ref{the: K-theory version implies KH-version} imply:

\begin{enumerate}

\item The $\KH$-Isomorphism Conjecture~\ref{con: KH-Isomorphism-Conjecture} is true for $G$
      and $R$;

\item The canonical map  $K_n(RG) \to \KH_n(RG)$ is bijective for $n \in \IZ$;

\item $N^p K_n(RG) = 0$ for $p \ge 1$ and $n \in \IZ$, see \cite[Proposition~7.4]{Bartels-Farrell-Jones-Reich(2004)}.

\end{enumerate}

\end{remark}

\begin{remark}[Injectivitiy of the $\KH$-assembly map]
In many cases injectivitiy of the assembly map
\[
H_n^G(\EGF{G}{\calfin};\bfK_R) \to K_n(RG)
\]
is proven by construction of a spectrum $\bfT{(R,G)}$ and a map of spectra $\bfK^{-\infty}{RG} \to \bfT{(R,G)}$
such that for many groups the composition of the assembly map on the level of spectra
with this map is a weak equivalence.
The construction of $\bfK^{-\infty}{RG} \to \bfT{(R,G)}$ is always natural in the coefficient ring $R$.
Therefore applying the arguments of the proof of
Theorem~\ref{the: K-theory version implies KH-version}~\ref{the: K-theory version implies KH-version (polynom)}
we can use $\bfT{(R[\bullet],G)}$ to split the $\KH$-assembly map in this cases.
This proves that the $\KH$-assembly map is split injective for groups $G$ of finite asymptotic dimension that admit
a finite model for $BG$ \cite{Bartels(2003a)} and for groups $G$ for which $\EGF{G}{\calfin}$ has a
compactification with certain properties \cite{Rosenthal(2003)}.
\end{remark}


\typeout{----------- Non-connective Waldhausen nil-spectra --------------}
\section{Non-connective Waldhausen Nil}
\label{sec: Non-connective Waldhausen Nil}

Before we can show that $H^?_*(-;\bfKH_R)$ has the tree property,
we will need to recall Waldhausen's work on $K$-theory of generalized free
products \cite{Waldhausen(1978a)}. We start with Waldhausen's $\Nil$-groups.

\begin{definition}[$\Nil$-categories]
Let $R$ be a ring and $X,Y,Z,W$ be  $R$-bimodules.
\begin{enumerate}
\item
The category $\NIL(R;X,Y)$ has
as objects quadruples $(P,Q,p,q)$, where $P$ and $Q$ are finitely generated
projective $R$-modules and $p\colon P \to Q \ox_R X$, $q: Q \to P \ox_R Y$ are
$R$-linear maps subject to the following {\em nilpotence} condition:
Let $P_0 = Q_0 = 0$, $P_{n+1} = p^{-1}(Q_n \ox_R X)$ and
$Q_ {n+1} = q^{-1}(P_n \ox_R Y)$. It is required that for
sufficient large $N$,
$P = P_N$ and $Q = Q_N$.
\item
The category $\NIL(R;X,Y,Z,W)$ has as objects
quadruples $(P,Q,p,q)$, where $P$ and $Q$ are finitely
generated projective $R$-modules
and $p\colon P \to Q \ox_R X \oplus P \ox_R Z$,
$q: Q \to P \ox_R Y \oplus Q \ox_R W$
are $R$-linear maps subject to the following
{\em nilpotence} condition: Let $P_0 = Q_0 = 0$,
$P_{n+1} = p^{-1}(Q_n \ox_R X \oplus P_n \ox_R Z)$ and
$Q_{n+1} = q^{-1}(P_n \ox_R Y \oplus Q_n \ox_R W)$.
It is required that for sufficient large $N$, $P = P_N$ and
$Q = Q_N$.
\end{enumerate}
Morphisms are in both cases $R$-linear maps
$P \to P'$, $Q \to Q'$ that are compatible with $p$, $p'$, $q$ and $q'$.
Both  categories are exact categories, where sequences are exact
whenever they map to exact sequences of modules under
$(P,Q,p,q) \mapsto P$ and $(P,Q,p,q) \mapsto Q$.
\end{definition}

\begin{remark}
\label{rem: functoriality-Nil}
Let $f_R\colon R \to S$ be a map of rings and
$f_X\colon X \to X'$, $f_Y\colon Y \to Y'$
be maps over $f_R$, i.e.\ $X'$ and $Y'$ are $S$-bimodules,
$f_X(rxr') = f_R(r) f_X(x) f_R(r')$
and similar for $f_Y$. Then $(f_R,f_X,f_Y)$ induce an exact  functor
$\NIL(R;X,Y) \to \NIL(S;X',Y')$ sending $(P,Q,p,q)$ to
$(P \ox_R S, Q \ox_R S, p_S, q_S)$ where $p_S$ and $q_S$ are the
canonical maps.
For example, $p_S$ is the composition
\begin{eqnarray*}
P \ox_R S \to Q \ox_R X \ox_R S \to Q \ox_R X' \cong Q \ox_R S \ox_S X',
\end{eqnarray*}
where the first map uses $p$ and the second uses $f_X$ and left multiplication
of $S$. In particular, we get  a functor
$S \ox -\colon \NIL(R;X,Y) \to \NIL(S \ox R; S \ox X, S \ox Y)$.
If $f\colon S \to S'$ is a map of rings, then we get another functor
$f_*\colon \NIL(S \ox R; S \ox X, S \ox Y) \to
           \NIL(S' \ox R; S' \ox X, S' \ox Y)$.
The functoriality of $\NIL(R;X,Y,Z,W)$ is similar.
\end{remark}

We review next \cite[2.5]{Wagoner(1972)} in a slightly more modern
language and discuss applications to Waldhausen Nil-categories.

A \emph{sum ring} is a ring $S$ together with elements $v$, $\bar{v}$,
$u$ and $\bar{u}$
of $S$ such that $u \bar{u} = 1$, $v \bar{v} = 1$ and
$\bar{v} v + \bar{u} u = 1$. This implies that $u \bar{v} = 0$ and
$v \bar{u} = 0$. Moreover, the map $f_\oplus\colon S \oplus S \to S$
defined by $(r,s) \mapsto \bar{u} r u + \bar{v} s v$ is a ring homomorphism.
Let $M,N$ be $S$-modules. Denote by $(M,N)$ the direct sum $M \oplus N$
considered
as an $S \oplus S$-module. The $S$-modules $(M,N) \ox_{f_\oplus} S$ and
$M \oplus N$ (considered as an $S$-module as usual) are naturally isomorphic.
Such an isomorphism and its inverse are given by
\[
\begin{array}{rcccccl}
M \oplus N & \ni & m \oplus n & \mapsto & (m,n) \ox (\bar{u} + \bar{v}) &
                                              \in & (M,N) \ox_{f_\oplus} S \\
(M,N) \ox_{f_\oplus} S & \ni & (m,n) \ox r & \mapsto & m u r \oplus  n v r &
                                              \in & M \oplus N.
\end{array}
\]
An \emph{infinite sum ring} is a sum ring together with a ring endomorphism $f_\infty$
such that $f_\oplus(r,f_\infty(r)) = f_\infty(r)$.

\begin{remark} \label{rem: swindle-on-infinite-sum}
The functor $M \mapsto M \ox_{f_\infty} S$ is an Eilenberg swindle
on the  category $\calp_S$ of finitely
generated projective modules over such an infinite sum ring. Indeed,
\begin{eqnarray*}
M \ox_{f_\infty} S & = & M \ox_{f_\oplus \circ (\id_S,f_\infty)} S \\
     & \cong & M \ox_{(\id_S,f_\infty)} (S \oplus S) \ox_{f_\oplus} S \\
     & \cong & ( M , M \ox_{f_\infty} S) \ox_{f_\oplus} S \\
     & \cong & M \oplus (M \ox_{f_\infty} S).
\end{eqnarray*}
The same swindle applies to
Waldhausen's Nil-categories:
Fix an infinite sum ring $S$. Let $X$ and $Y$ be
bimodules over another  ring $R$.
Then the endofunctor $(f_\infty)_*$ is equal to the composition
\[
\xymatrix
{
\NIL(S \ox R; S \ox X, S \ox Y) \ar[d]_{(\id,f_\infty)_*} \\
\NIL((S \oplus S) \ox R; (S \oplus S) \ox X, (S \oplus S) \ox Y)
      \ar[d]_{(f_\oplus)_*} \\
\NIL(S \ox R; S \ox X, S \ox Y) \\
}
\]
Using the natural isomorphism from above there is a natural transformation
from this composition to $\id \oplus (f_\infty)_*$. Thus, $(f_\infty)_*$
is an Eilenberg swindle.
This swindle is compatible with the two forgetful functors
$\NIL(S \ox R; S \ox X, S \ox Y) \to \calp_{S \ox R}$.
Analogous considerations apply to
$\NIL(S \ox R; S \ox X, S \ox Y, S \ox Z, S \ox Z, S \ox W)$.
\end{remark}

The {\em cone ring} $\conering \IZ$ of $\IZ$ is the ring of
column and row finite
$\IN \x \IN$-matrices over $\IZ$, i.e.\ matrices such that every column and
every row contains only finitely
many non-zero entries. The {\em suspension ring} $\Sigma \IZ$ is the quotient
of $\conering \IZ$
by the ideal of finite matrices. For an arbitrary ring $R$ we define
$\conering R = \conering \IZ \ox R$ and $\Sigma R = \Sigma \IZ \ox R$. We will
view $\conering$ and $\Sigma$ as functors.
Every bijection $\IN \to \IN \x \IN$ induces the structure of an infinite
sum ring on the cone ring $\conering R$, cf.~\cite[p.355]{Wagoner(1972)}.
We can consider $\conering R$ as a
subring of the ring of column and row finite
$\IN \x \IN$-matrices over $R$. However, this inclusion is not always an
equality.

Next we want to define a non-connective spectrum associated to
Waldhausen's Nil-categories.
First we recall the construction for $K$-theory.
Denote by $\IK R$ the $K$-theory space of a ring
(obtained for example by applying
Waldhausen's $S_\bullet$-construction to the category $\calp_R$ of finitely
generated projective $R$-modules).
The $n$-th space of the spectrum $\bfK^{-\infty} R$ is by definition $\IK \Sigma^n R $.
The composition  $\IK R \to \IK \conering R \to \IK \Sigma R$ is  constant.
The choice of an bijection $\IN \to \IN \x \IN$ gives an Eilenberg swindle on
$\conering R$, cf.\ Remark~\ref{rem: swindle-on-infinite-sum}.
If we fix such a bijection we get
a functorial way of contracting  $\IK \conering R$ to the basepoint.
This induces
the structure maps $\Sigma (\IK \Sigma^n R) \to \IK \Sigma^{n+1} R$.

\begin{definition}[Non-connective $\Nil$-spectra] \label{def: non-connective-Nil}
Let $R$ be a ring and $X$ and $Y$ be $R$-bimodules. The
(non-connective) spectrum $\bfNil^{-\infty}(R;X,Y)$ has
$\IK \NIL(\Sigma^n R; \Sigma^n X, \Sigma^n Y)$ as its
nth space.
Here $\Sigma^n X = \Sigma^n \IZ \ox X$ and we define similarly
$\Sigma^n Y$, $\conering X$ and $\conering Y$.
The structure maps are defined in an
analogous ways as for the non-connective $K$-theory spectrum:
The functoriality
discussed in Remark~\ref{rem: functoriality-Nil} allows us to consider the
(constant) composition
\[
\NIL(C;A',B') \to \NIL(\conering C; \conering A', \conering B')
     \to \NIL(\Sigma C; \Sigma A', \Sigma B').
\]
The structure maps for $\bfNil^{-\infty}(R;X,Y)$ are now defined using the Eilenberg
swindle on the second category discussed in
Remark~\ref{rem: swindle-on-infinite-sum}.
Similarly, we define a (non-connective) spectrum
$\bfNil^{-\infty}(R;X,Y,Z,W)$ with
\[
\IK \NIL(\Sigma^n R; \Sigma^n X, \Sigma^n Y, \Sigma^n Z, \Sigma^n W)
\]
as its nth space.
\end{definition}

An inclusion $\alpha\colon C \to A$ of rings is called {\em pure} if
$A = \alpha(C) \oplus A'$
as $C$-bimodules. It is called {\em pure and free} if in addition $A'$ is
free as a left $C$-module. If $H \to G$ is an inclusion of groups,
then the inclusion $RH \to RG$ of rings is
pure and free.
The following observation is straight forward.

\begin{lemma}
If $\alpha$ is pure (and free) then $\Sigma \alpha$ and
$\conering \alpha$ are also
pure (and free).
\end{lemma}

Let $\alpha\colon C \to A$ and $\beta\colon C \to B$ be both pure.
The ring $R = A *_C B$,
the {\em free product of $A$ and $B$, amalgamated at $C$
(w.r.t.\ $\alpha,\beta$)}, is defined by the push-out
\[
\xymatrix
{
C \ar[r]^\alpha \ar[d]_\beta & A \ar[d] \\
B \ar[r] & R.
}
\]
For group rings this corresponds to amalgamated products of groups.

\begin{lemma} \label{cone-commutes-w-amalg}
The cone respectively suspension ring of
$A *_C B$
is naturally isomorphic to
$\conering A *_{\conering C} \conering B$ respectively
$\Sigma A *_{\Sigma C} \Sigma B$.
\end{lemma}

\begin{proof}
This follows from the universal property.
\end{proof}

Let $\alpha, \beta\colon C \to A$ be pure and free. The
{\em Laurent extension w.r.t.\ $\alpha$ and $\beta$} is the universal ring
$R = A_{\alpha,\beta} \{ t^{\pm 1} \}$ that contains $A$ and an invertible
element $t$ and satisfies
\[
\alpha(c) t = t \beta(c) \; \mbox{for} \; c \in C.
\]
Existence is explained in \cite[p.149]{Waldhausen(1978a)}.
For group rings this corresponds to HNN-extensions.

\begin{lemma} \label{cone-commutes-w-Laurent}
The cone respectively  suspension ring of
$A_{\alpha,\beta} \{ t^{\pm 1} \}$
is naturally isomorphic to
$\conering A_{\conering \alpha,\conering \beta} \{ T^{\pm 1} \}$ respectively
$\Sigma A_{\Sigma \alpha,\Sigma \beta} \{ T^{\pm 1} \}$.
\end{lemma}

\begin{proof}
This follows from the universal property.
\end{proof}


\typeout{----------- Waldhausen's cartesian squares ----------}
\section{Waldhausen's cartesian squares}
\label{subsec: Waldhausen's cartesian squares}

Let $\alpha\colon C \to A$ and $\beta\colon C \to B$ be pure and free.
Write $A  = \alpha(C) \oplus A'$ and  $B  = \beta(C) \oplus B'$ as
$C$-bimodules. Let $R = A *_C B$.
Consider the square
\begin{eqnarray}
&
\comsquare{\NIL(C;A',B')}{}{\calp_A \times \calp_B}{}{}{\calp_C}{}{\calp_R^*}
&
\label{free-diagram}
\end{eqnarray}
The two functors starting at the upper left hand corner are defined by
sending $(P,Q,p,q)$
to $(P \oplus Q) \in \calp_C$ respectively to $(P \ox_\alpha A, Q \ox_\beta B)$.
The category
$\calp^*_R$ is defined in \cite[p.205]{Waldhausen(1978a)}. It is a cofinal full
subcategory of $\calp_R$ and contains all finitely generated free modules.
There is an obvious natural transformation between the two ways to go through
the diagram. However, there is also a not quite so obvious more complicated
natural transformation that makes use of $p$ and $q$,
cf.~\cite[1.4,11.3]{Waldhausen(1978a)}:
Let
$i_P\colon P \ox_C A \ox_A R \to P \ox_C R$ and
$i_Q\colon Q \ox_C B \ox_B  R \to Q \ox_C R$ be the
natural isomorphisms. Define $N$ by the commutative diagram
$$\comsquare{P}{p}{Q \otimes_C A}{}{}
{P \otimes_C A \otimes_A R}{N}{Q \otimes_C R}$$
and similarly $M\colon Q \ox_C B \ox_B R \to P \ox_C R$. Then
 \begin{eqnarray*}
 \left(
 \begin{array}{cc}
 i_P & M \\ N & i_Q
 \end{array}
 \right)
 \end{eqnarray*}
is an isomorphism and defines the more complicated natural transformation.
It is a result of Waldhausen \cite[11.3]{Waldhausen(1978a)},
that applying $\IK$ to (\ref{free-diagram}) yields a homotopy cartesian square
(w.r.t.\ the homotopy induced by the more complicated natural
transformation).
We will need a non-connective version of Waldhausen's result.

\begin{theorem}{\bf (Non-connective versions of Waldhausen's homotopy cartesian squares for
    amalgamation)}
\label{the: free-seq}
We have the following diagram of spectra
$$\comsquare{\bfNil^{-\infty}(C;A',B')}{}{\bfK^{-\infty} A \wedge \bfK^{-\infty} B }{}{}
{\bfK^{-\infty} C}{}{\bfK^{-\infty} R}
$$
The more complicated natural
transformations combine to a homotopy between the two
ways to go through this diagram. The diagram is homotopy cartesian w.r.t.\ this
homotopy.
\end{theorem}

\begin{proof}
The diagram of spectra is obtained from (\ref{free-diagram}) by tensoring
everything in sight by $\Sigma^n \IZ$ (and applying $\IK$). We need to check
compatibility with the structure maps. Those come from an intermediate diagram
where we apply $\conering \IZ \ox -$ and use an Eilenberg swindle.
This Eilenberg
swindle happens on the left of this tensor product,
while everything else happens on the right. This proves compatibility with the
structure maps. If we use $\calp^*_R$ rather then $\calp_R$ then the diagram
is homotopy cartesian by Waldhausen's result and
Lemma~\ref{cone-commutes-w-amalg}. However, since the former
category contains all finitely generated
free modules and we use non-connective $K$-theory we can also use $\calp_R$.
\end{proof}

\begin{remark}{\bf (Waldhausen Nil for amalgamations vanishes for regular coherent rings)}
\label{rem: nil vanishes(amalgamation)}
Waldhausen proved that for a  regular coherent ring $C$,
the functor $\calp_C \x \calp_C \to \NIL(C;A',B')$ defined by $(P,Q) \mapsto (P,Q,0,0)$
induces an isomorphism in connective $K$-theory, \cite[12.2]{Waldhausen(1978a)}.
A priori, this does not immediately imply that the induced map
$\alpha :
\bfK^{-\infty} C \vee \bfK^{-\infty} C \to \bfNil^{-\infty}(C;A',B')$
is a weak equivalence
because it is not clear whether $\Sigma C$ is again regular coherent.
However, if $C$ is regular or
more generally, if $C$ is a group ring with a regular coefficient ring over a regular coherent group,
then $\alpha$ is a weak equivalence.
This can be seen as follows.
The functor $(P,Q,p,q) \mapsto (P,Q)$ splits $\alpha$,
thus $\alpha$ will injective on homotopy groups.
To prove surjectivity we use the fact that there is an in $C$ natural map of rings $C[\IZ^n] \to \Sigma^n C$,
that is naturally split surjective in connective $K$-theory,
\cite[Section 6]{Wagoner(1972)}.
We get the following commutative diagram
\begin{equation}
\label{eq: laurent to suspension}
\comsquare{\calp_{C[\IZ^n]} \x \calp_{C[\IZ^n]}}{}{\NIL(C[\IZ^n];A'[\IZ^n],B'[\IZ^n])}
          {}{}
          {\calp_{\Sigma^n C} \x \calp_{\Sigma^n C}}{}{\NIL(\Sigma^n C;\Sigma^n A',\Sigma^n B')}
\end{equation}
Apply $C[\IZ^n] \to \Sigma^n C$ map to the long exact sequence obtained from \eqref{free-diagram}
by \cite[11.3]{Waldhausen(1978a)}.
A little diagram chase in the resulting ladder diagram shows that the right
vertical map in \eqref{eq: laurent to suspension}
is surjective in connective $K$-theory.
The assumptions on $C$ imply that $C[\IZ^n]$ is regular coherent.
Therefore the top horizontal map in \eqref{eq: laurent to suspension}
is an isomorphism in connective $K$-theory.
Therefore the bottom horizontal map in \eqref{eq: laurent to suspension} is
also surjective in connective $K$-theory.
This implies that $\alpha$ is surjective on homotopy groups.
\end{remark}

Next we discuss the analogous cartesian square for Laurent extensions.
Let $\alpha, \beta\colon C \to A$ be pure and free and
$R = A_{\alpha,\beta} \{ t^{\pm 1} \}$.
We denote by $\iota\colon A \to R$ the inclusion.
Write  $A = \alpha(C) \oplus A'$ and $A = \beta(C) \oplus A''$ as
$C$-bimodules.
Consider the square
\begin{eqnarray}
& \comsquare{\NIL(C; _\alpha A'_\alpha, _\beta A''_\beta, _\beta
  A_\alpha, _\alpha A_\beta)}{}
{\calp_A}{}{\iota_*}
{\calp_C}{(\iota \circ \alpha)_*}{\calp^*_R}
&
\label{L-diagram}
\end{eqnarray}
Here we we use $\alpha$ and $\beta$ to indicate the
$C$-bimodule structures.
The two functors starting at the upper right hand corner
are defined by sending $(P,Q,p,q)$ to $(P \oplus Q) \in \calp_C$ respectively
to $(P \ox_\alpha A \oplus Q \ox_\beta A)$. The category $\calp^*_R$ is defined
in \cite[p.205]{Waldhausen(1978a)}. It is a cofinal full
subcategory of $\calp_R$ and contains all finitely generated free modules.
As before
there is an obvious and a more complicated  natural transformation
between the two ways to go through the
diagram \cite[2.4,12.3]{Waldhausen(1978a)}:
Let $i_P\colon P \ox_\alpha A \ox_A \ox R \to P \ox_\alpha R$ and
$i_Q\colon Q \ox_\beta A \ox R \to Q \ox_\alpha R$ denote the canonical
isomorphisms.
(Here $i_Q$ uses an extra $t$, i.e.\
$i_Q(y \ox a \ox r) = y \ox t a r$.)
These isomorphisms give the obvious natural transformation.
The more complicated natural transformation is obtained by adding a
nilpotent term which we review next. Write $p = p_0 + p_1$, where
$p_0\colon P \to P \ox _\beta A_\alpha$ and
$p_1\colon P \to Q \ox _\alpha A'_\alpha$. Define $N_0$ and $N_1$ by
the commutative diagrams
$$\comsquare{P}{p_0}{P \otimes_\beta A_\alpha}{}{}
{P \otimes_\alpha A \otimes_A R}{N_0}{P \otimes_\alpha R}
\hspace{12mm}
\comsquare{P}{p_1}{Q \otimes_\beta A'_\alpha}{}{}
{P \otimes_\alpha A \otimes_A R}{N_1}{Q \otimes_\alpha R}
$$
(The second vertical arrow is $x \ox a \mapsto x \ox t a$.)
Write $q = q_0 + q_1$, where
$q_0\colon Q \to Q \ox _\alpha A_\beta$ and
$q_1\colon Q \to P \ox _\beta A''_\beta$. Define $M_0$ and $M_1$ by
the commutative diagrams
$$\comsquare{Q}{q_0}{Q \otimes_\alpha A_\beta}{}{}
{Q \otimes_\beta A \otimes_A R}{M_0}{Q \otimes_\alpha R}
\hspace{12mm}
\comsquare{Q}{q_1}{P \otimes_\alpha A''_\beta}{}{}
{Q \otimes_\beta A \otimes_A R}{M_1}{P \otimes_\alpha R}
$$
(The forth vertical arrow is $x \ox a \mapsto x \ox t a$.)
The more complicated natural transformation is then given by
the isomorphism
\[
\left(
\begin{array}{cc}
i_P & 0 \\ 0 & i_Q
\end{array}
\right)
+
\left(
\begin{array}{cc}
N_0 & M_1 \\ N_1 & M_0
\end{array}
\right).
\]
It is a result of
Waldhausen \cite[12.3]{Waldhausen(1978a)},
that applying $\IK$ to (\ref{L-diagram}) yields a homotopy cartesian square
(w.r.t.\ the homotopy induced by
the more complicated natural transformation).
The arguments used to prove Theorem~\ref{the: free-seq} can be
used to prove a non-connective version of this result.

\begin{theorem}{\bf (Non-connective versions of Waldhausen's homotopy cartesian squares for
Laurent extensions)}  \label{the: Laurent-square}
We have the following diagram of spectra
$$\comsquare{\bfNil^{-\infty}(C; _\alpha A'_\alpha, _\beta A''_\beta, _\beta A_\alpha,
                                                _\alpha A_\beta)}{}{\bfK^{-\infty}(A)}
{}{\bfK^{-\infty}(\iota)}
{\bfK^{-\infty}(C) }{\bfK(\iota \circ \alpha)}{\bfK^{-\infty}(R)}$$
The more complicated natural transformations combine to
a homotopy between the two ways to go through this diagram. The diagram
is homotopy cartesian w.r.t.\ this homotopy.
\end{theorem}

\begin{remark}{(\bf Waldhausen Nil for Laurent extensions vanishes for regular coherent rings)}
\label{rem: nil vanishes(Laurent)}
The reasoning in Remark~\ref{rem: nil vanishes(amalgamation)} also applies to
$\Nil(C; _\alpha A'_\alpha, _\beta A''_\beta, _\beta A_\alpha, _\alpha A_\beta)$.
If $C$ is a group ring with a regular coefficient ring over a regular coherent group,
then the functor $(P,Q) \mapsto (P,Q,0,0)$ induces a weak equivalence
$ \bfK^{-\infty} C \vee \bfK^{-\infty} C \to
  \bfNil^{-\infty}(C; _\alpha A'_\alpha, _\beta A''_\beta, _\beta A_\alpha, _\alpha A_\beta)$.
\end{remark}

%


\typeout{----------- The tree property for KH ----------}
\section{The tree property for Homotopy K-Theory}
\label{sec: The tree property for KH}

This section contains the proof of the following result.

\begin{theorem}[Continuity and tree-property for $H^?_*(-;\bfKH_R)$]
\label{the: tree-property-for-KH}
The equivariant homology theory $H^?_*(-;\bfKH_R)$
is continuous and has the tree property.
\end{theorem}

Let $X$, $Y$,$Z$ and $W$ be bimodules over $R$.
Consider the simplicial spectra
$\ULn \mapsto \bfNil(\IZ[\ULn] \ox R;\IZ[\ULn] \ox X ,\IZ[\ULn] \ox Y)$
and
$\ULn \mapsto \bfNil(\IZ[\ULn] \ox R;\IZ[\ULn] \ox X, \IZ[\ULn] \ox Y,
                                     \IZ[\ULn] \ox Z, \IZ[\ULn] \ox W)$.
We will denote the realization of these simplicial spectra by $\bfNH(R;X,Y)$
or $\bfNH(R;X,Y,Z,W)$ respectively.
However, the point here is that this process kills the additional
information in Waldhausen's $\Nil$-groups.

\begin{proposition} \label{pro: no-nil}
$ $
\begin{enumerate}
\item \label{no-nil:amalg}
The functor $\calp_R \x \calp_R \to \NIL(R;X,Y)$ defined by
$(P,Q) \mapsto (P,Q,0,0)$ induces an equivalence
$\bfKH (R) \vee \bfKH(R) \to \bfNH(R;X,Y)$.
\item \label{no-nil:Laurent}
The functor $\calp_R \x \calp_R \to \NIL(R;X,Y,Z,W)$ defined by
$(P,Q) \mapsto (P,Q,0,0)$ induces an equivalence
$\bfKH (R) \vee \bfKH(R) \to \bfNH(R;X,Y,Z,W)$.
\end{enumerate}
\end{proposition}

\begin{proof}
We prove only \ref{no-nil:amalg}, the proof of \ref{no-nil:Laurent} is similar.
It suffices to show that the functor $(P,Q,p,q) \mapsto (P,Q,0,0)$ mapping the
simplicial category
$\categoryN(\bullet) = \NIL(\IZ[\bullet] \ox R;
           \IZ[\bullet] \ox X ,\IZ[\bullet] \ox Y)$
to itself is simplicially  homotopic to the identity. Such a homotopy
$\categoryN(\bullet) \x  \Delta^1_\bullet \to \categoryN(\bullet)$
is given by
\[
(P,Q,p,q) \mapsto (P, Q, \sum_{j \in f^{-1}(0)} t_j \ox p,
                         \sum_{j \in f^{-1}(0)} t_j \ox q)
\]
where $f\colon \ULn \mapsto \UL{1}$.
\end{proof}

The homotopy algebraic $K$-theory of a free product or a Laurent extension
does therefore not involve $\Nil$-groups.

\begin{theorem}{\bf (Homotopy cartesian squares for Homotopy $K$-theory)} \label{cartesian-square-for-homotopy-K}
\begin{enumerate}
\item \label{cartesian-square-for-homotopy:free}
Consider the free product $R = A *_C B$ (w.r.t.\ pure
and free maps $\alpha\colon C \to A$ and $\beta\colon C \to B$).
Then the commutative diagram
\begin{eqnarray*}
\comsquare
{\bfKH(C) \vee \bfKH(C)} {\bfKH(\alpha) \vee \bfKH(\beta)} {\bfKH(A) \vee \bfKH(B)}
{\bfid \vee \bfid} {\bfKH(\iota_A) \vee \bfKH(\iota_B)}
{\bfKH(C)} {\bfKH(\alpha \circ \iota_A)} {\bfKH(R)}
\end{eqnarray*}
is homotopy cartesian. Here $\iota_A$ and $\iota_B$ are the
obvious inclusions of rings.
\item \label{cartesian-square-for-homotopy:Laurent}
Consider the Laurent extension
$R = A_{\alpha, \beta} \{ t^{\pm 1} \}$ (w.r.t.\ pure
and free maps $\alpha\colon C \to A$ and $\beta\colon C \to A$).
Then the diagram
\begin{eqnarray*}
\comsquare
{\bfKH(C) \vee \bfKH(C)} {\bfKH(\alpha) \vee \bfKH(\beta)} {\bfKH(A)}
{\bfid \vee \bfid} {\bfKH(\iota)}
{\bfKH(C)} {\bfKH(\iota \circ \alpha)}  {\bfKH(R)}
\end{eqnarray*}
is homotopy cartesian w.r.t.\ the obvious homotopy between the two
ways to go through this diagram,
cf. Section~\ref{subsec: Waldhausen's cartesian squares}.
Here $\iota$ is  the obvious inclusions of rings.
\end{enumerate}
\end{theorem}

\begin{proof}
The realization of a simplicial diagram {\em of spectra} that is degreewise
homotopy cartesian is again homotopy cartesian.
Thus the result follows by
combining Theorems~\ref{the: free-seq} or \ref{the: Laurent-square} respectively and
Proposition~\ref{pro: no-nil}.
There is no longer a difference between the obvious
and the more complicated natural transformation,
since we got rid of the nil-categories.
\end{proof}

\begin{proof}[Proof of Theorem~\ref{the: tree-property-for-KH}]
We use Lemma~\ref{lem: tree property for calh^?_*(-,bfE)}.
We discuss first continuity, i.e. condition
\ref{lem: tree property for calh^?_*(-,bfE)}
\ref{lem: tree property for calh^?_*(-,bfE): continuity}.
For $\bfK_R$ this follows from
the compatibility of $K$-theory with directed colimits
\cite{Quillen(1973)}.
Since realizations of simplicial spectra  commutes with directed
colimits, this implies condition
\ref{lem: tree property for calh^?_*(-,bfE)}
\ref{lem: tree property for calh^?_*(-,bfE): continuity}.

Next we discuss the tree property, i.e. conditions
\ref{lem: tree property for calh^?_*(-,bfE)}
\ref{lem: tree property for calh^?_*(-,bfE): amalgamated product}
and \ref{lem: tree property for calh^?_*(-,bfE): HNN}.
Note that $\bfKH_R$ is a $2$-functor as discussed in
Remark~\ref{rem: smaller diagrams}. For $\bfK_R$ this holds
since natural equivalences between functors $F$ and $G$ induce naturally
a homotopy from $\bfK^{-\infty}(F)$ to $\bfK^{-\infty}(G)$. This homotopy is preserved under
realization.
Now observe that the obvious homotopy in
Theorem~\ref{cartesian-square-for-homotopy-K}
\ref{cartesian-square-for-homotopy:Laurent}
is in the case of an HNN-extension the homotopy coming from
conjugation as discussed in Remark~\ref{rem: smaller diagrams}
\ref{rem: smaller diagrams:HNN}.
Thus conditions \ref{lem: tree property for calh^?_*(-,bfE)}
\ref{lem: tree property for calh^?_*(-,bfE): amalgamated product}
and \ref{lem: tree property for calh^?_*(-,bfE): HNN}
are satisfied for $\bfKH_R$ by
Theorem~\ref{cartesian-square-for-homotopy-K}
and Remark~\ref{rem: smaller diagrams}.
\end{proof}

Using Theorems~\ref{the: free-seq} and \ref{the: Laurent-square} and
Remarks~\ref{rem: nil vanishes(amalgamation)} and \ref{rem: nil vanishes(Laurent)}
the above arguments also prove a version of Theorem~\ref{the: tree-property-for-KH}
for algebraic $K$-theory.

\begin{theorem}[Continuity and tree-property for $H^?_*(-;\bfK_R)$]
\label{the: regular tree property for K}
The equivariant homology theory $H^?_*(-;\bfK_R)$
is continuous and if $R$ is regular then it has the regular tree property.
\end{theorem}

Now we can  finish the proof of the various results stated in the Introduction.
We start with Theorem~\ref{the: inhertitance for KH} and Theorem~\ref{the: inhertiance for Farrell Jones}.
The property (FIN) respectively (VCYC) hold for trivial reasons.
Similar (SUB) in the Fibered case is a formal consequence of the Definitions,
compare Lemma~\ref{lem: basic inheritance property of fibered conjecture}.
The property (COL) is a consequence of Proposition~\ref{prop: isomorphism conjecture is stable under colim},
Theorem~\ref{the: tree-property-for-KH} and
Theorem~\ref{the: regular tree property for K}.
The properties (TREE) respectively (TREE$_\calr$) follow from
Theorem~\ref{the: Isomorphisms Conjecture and actions on trees} and
Theorem~\ref{the: tree-property-for-KH} respectively Theorem~\ref{the: regular tree property for K}.
Theorem~\ref{the: extensions and actions with finite stabilizers} is a consequence of
Corollary~\ref{cor: extensions by groups acting on trees with finite stabilizers}
and Theorem~\ref{the: tree-property-for-KH}.
Theorem~\ref{the: conclusions for calc and K-theoretic Farrell-Jones} is a consequence of
Theorem~\ref{the: inhertitance for KH}, Proposition~\ref{pro: comparing-K-assembly-to-KH-assembly}
and Remark~\ref{rem: assumption of propositions hold in cases}.
It remains to prove Proposition~\ref{pro: groups belonging to calc}.
\begin{proof}
\ref{pro: groups belonging to calc: one relator groups} This result is stated in
\cite[page 133]{Oyono-Oyono(2001b)}. We give an outline of the proof.
We consider first a group $G$ which possesses
a finite presentation with one relation.
Let $r$ be the number of generators appearing in the word describing the
relation. If $r \le 1$, then $G$ is the amalgamated product of
a free group and a finite
cyclic group. Obviously any finite group and $\IZ$ belong to $\calc_0$ and $\calc_0$ is closed under free products.
Hence $G$ belongs to $\calc_0$.
It remains to treat $r \ge 2$.
Here we use induction over the length $l$ of the word describing the relation.
In our case $l \geq 2$.
Then $G$ acts on a tree with stabilizers which are subgroups
of one-relator
groups whose relation have length $\le (l-1)$,
\cite[Theorem 7.7]{Bieri(1976a)} and hence belong to $\calc_0$ by the induction hypothesis.
Therefore $G$ belongs to $\calc_0$.
For a general one-relator group $G$ there are finitely generated
subgroups $G_i$ which are free or one-relator groups, such that
$G$ is is the directed colimit over the $G_i$.
\\[2mm]
\ref{pro: groups belonging to calc: poly-free}
Let $\{1\} = G_0 \subseteq G_1 \subseteq \ldots \subseteq G_n = G$ be a sequence
of subgroups such that $G_{i-1}$ is normal in $G_i$ and the quotient
$G_i/G_{i-1}$ is free for $i = 1,2 \ldots , n$. We prove by induction
over $n$ that $G$ belongs to $\calc_0$.
The induction beginning $n = 0$ is trivial because of property (FIN),
the induction step done as follows.
We can write $G_n/G_{n-1}$ as a directed union of its finitely generated subgroups.
Hence $G_n$ is the directed union of the preimages of the finitely generated
subgroups of $G_n/G_{n-1}$. Since any finitely generated subgroup of a
free group $F$ is a finitely generated free group, it suffices to treat the case, where $G_n/G_{n-1}$ is
finitely generated free by property (COL).  Since $G_n/G_{n-1}$ acts
on a tree with trivial stabilizers, $G_n$ acts on a tree with
stabilizers
which are all isomorphic to $G_{n-1}$ and hence belong to $\calc_0$.
Hence $G$ belongs to $\calc_0$ by property (TREE).
\\[2mm]
\ref{pro: groups belonging to calc: 3-amifold groups},
\ref{pro: groups belonging to calc: surface groups}, and
\ref{pro: groups belonging to calc: submanifolds of S^3}
These follow from \cite[Theorem 17.5 on page 250]{Waldhausen(1978a)}.
\end{proof}

\def\cprime{$'$} \def\polhk#1{\setbox0=\hbox{#1}{\ooalign{\hidewidth
  \lower1.5ex\hbox{`}\hidewidth\crcr\unhbox0}}}


\begin{thebibliography}{10}

\bibitem{Bartels(2003b)}
A.~Bartels.
\newblock On the domain of the assembly map in algebraic {$K$}-theory.
\newblock {\em Algebr. Geom. Topol.}, 3:1037--1050 (electronic), 2003.

\bibitem{Bartels(2003a)}
A.~Bartels.
\newblock Squeezing and higher algebraic {$K$}-theory.
\newblock {\em $K$-Theory}, 28(1):19--37, 2003.

\bibitem{Bartels-Farrell-Jones-Reich(2004)}
A.~Bartels, T.~Farrell, L.~Jones, and H.~Reich.
\newblock On the isomorphism conjecture in algebraic {$K$}-theory.
\newblock {\em Topology}, 43(1):157--213, 2004.

\bibitem{Bartels-Lueck(2004)}
A.~Bartels and W.~L\"uck.
\newblock Induction theorems and isomorphism conjectures for {$K$}- and
  {$L$}-theory.
\newblock Preprintreihe SFB 478 --- Geometrische Strukturen in der Mathematik,
  Heft 331, M\"unster, arXiv:math.KT/0404486, 2004.

\bibitem{Bartels-Reich(2005)}
A.~Bartels and H.~Reich.
\newblock On the {F}arrell-{J}ones conjecture for higher algebraic
  {$K$}-theory.
\newblock {\em J. Amer. Math. Soc.}, 18(3):501--545 (electronic), 2005.

\bibitem{Bass(1968)}
H.~Bass.
\newblock {\em Algebraic ${K}$-theory}.
\newblock W. A. Benjamin, Inc., New York-Amsterdam, 1968.

\bibitem{Bass-Heller-Swan(1964)}
H.~Bass, A.~Heller, and R.~G. Swan.
\newblock The {W}hitehead group of a polynomial extension.
\newblock {\em Inst. Hautes \'Etudes Sci. Publ. Math.}, (22):61--79, 1964.

\bibitem{Bieri(1976a)}
R.~Bieri.
\newblock {\em Homological dimension of discrete groups}.
\newblock Mathematics Department, Queen Mary College, London, 1976.
\newblock Queen Mary College Mathematics Notes.

\bibitem{Cappell(1974c)}
S.~E. Cappell.
\newblock Unitary nilpotent groups and {H}ermitian ${K}$-theory. {I}.
\newblock {\em Bull. Amer. Math. Soc.}, 80:1117--1122, 1974.

\bibitem{Davis-Lueck(1998)}
J.~F. Davis and W.~L{\"u}ck.
\newblock Spaces over a category and assembly maps in isomorphism conjectures
  in ${K}$- and ${L}$-theory.
\newblock {\em $K$-Theory}, 15(3):201--252, 1998.

\bibitem{Dunwoody(1979)}
M.~J. Dunwoody.
\newblock Accessibility and groups of cohomological dimension one.
\newblock {\em Proc. London Math. Soc. (3)}, 38(2):193--215, 1979.

\bibitem{Farrell-Jones(1993a)}
F.~T. Farrell and L.~E. Jones.
\newblock Isomorphism conjectures in algebraic ${K}$-theory.
\newblock {\em J. Amer. Math. Soc.}, 6(2):249--297, 1993.

\bibitem{Gersten(1973)}
S.~M. Gersten.
\newblock Higher {$K$}-theory of rings.
\newblock In {\em Algebraic $K$-theory, I: Higher $K$-theories (Proc. Conf.
  Seattle Res. Center, Battelle Memorial Inst., 1972)}, pages 3--42. Lecture
  Notes in Math., Vol. 341. Springer, Berlin, 1973.

\bibitem{Grayson(1976)}
D.~Grayson.
\newblock Higher algebraic ${K}$-theory. {I}{I} (after {D}aniel {Q}uillen).
\newblock In {\em Algebraic $K$-theory (Proc. Conf., Northwestern Univ.,
  Evanston, Ill., 1976)}, pages 217--240. Lecture Notes in Math., Vol. 551.
  Springer-Verlag, Berlin, 1976.

\bibitem{Juan-Pineda-Prassidis(2004)}
D.~Juan-Pineda and S.~Prassidis.
\newblock A controlled approach to the isomorphism conjecture.
\newblock preprint, arXiv:math.KT/0404513, 2004.

\bibitem{Lueck(2002b)}
W.~L{\"u}ck.
\newblock Chern characters for proper equivariant homology theories and
  applications to ${K}$- and ${L}$-theory.
\newblock {\em J. Reine Angew. Math.}, 543:193--234, 2002.

\bibitem{Lueck(2004a)}
W.~L\"uck.
\newblock Survey on classifying spaces for families of subgroups.
\newblock Preprintreihe SFB 478 --- Geometrische Strukturen in der Mathematik,
  Heft 308, M\"unster, arXiv:math.GT/0312378 v1, 2004.

\bibitem{Lueck-Reich(2004g)}
W.~L{\"u}ck and H.~Reich.
\newblock The {B}aum-{C}onnes and the {F}arrell-{J}ones conjectures in ${K}$-
  and ${L}$-theory.
\newblock Preprintreihe SFB 478 --- Geometrische Strukturen in der Mathematik,
  Heft 324, M\"unster, arXiv:math.GT/0402405, to appear in the
  $K$-theory-handbook, 2004.

\bibitem{Oyono-Oyono(2001b)}
H.~Oyono-Oyono.
\newblock Baum-{C}onnes conjecture and group actions on trees.
\newblock {\em $K$-Theory}, 24(2):115--134, 2001.

\bibitem{Pimsner-Voiculescu(1982)}
M.~Pimsner and D.~Voiculescu.
\newblock ${K}$-groups of reduced crossed products by free groups.
\newblock {\em J. Operator Theory}, 8(1):131--156, 1982.

\bibitem{Pimsner(1986)}
M.~V. Pimsner.
\newblock ${K}{K}$-groups of crossed products by groups acting on trees.
\newblock {\em Invent. Math.}, 86(3):603--634, 1986.

\bibitem{Quillen(1973)}
D.~Quillen.
\newblock Higher algebraic ${K}$-theory. {I}.
\newblock In {\em Algebraic $K$-theory, I: Higher $K$-theories (Proc. Conf.,
  Battelle Memorial Inst., Seattle, Wash., 1972)}, pages 85--147. Lecture Notes
  in Math., Vol. 341. Springer-Verlag, Berlin, 1973.

\bibitem{Rosenthal(2003)}
D.~Rosenthal.
\newblock Splitting with continuous control in algebraic {$K$}-theory.
\newblock {\em $K$-Theory}, 32(2):139--166, 2004.

\bibitem{Roushon(2004)}
S.~K. Roushon.
\newblock The {F}arrell-{J}ones isomorphism conjecture for 3-manifold groups.
\newblock preprint, arXiv:math.KT/0405211, 2004.

\bibitem{SauerJ(2002)}
J.~Sauer.
\newblock ${K}$-theory for proper smooth actions of totally disconnected
  groups.
\newblock Ph.D. thesis, 2002.

\bibitem{Serre(1980)}
J.-P. Serre.
\newblock {\em Trees}.
\newblock Springer-Verlag, Berlin, 1980.
\newblock Translated from the French by J.~Stillwell.

\bibitem{Switzer(1975)}
R.~M. Switzer.
\newblock {\em Algebraic topology---homotopy and homology}.
\newblock Springer-Verlag, New York, 1975.
\newblock Die Grundlehren der mathematischen Wissenschaften, Band 212.

\bibitem{Wagoner(1972)}
J.~B. Wagoner.
\newblock Delooping classifying spaces in algebraic ${K}$-theory.
\newblock {\em Topology}, 11:349--370, 1972.

\bibitem{Waldhausen(1978a)}
F.~Waldhausen.
\newblock Algebraic ${K}$-theory of generalized free products. {I}, {I}{I}.
\newblock {\em Ann. of Math. (2)}, 108(1):135--204, 1978.

\bibitem{Weibel(1981)}
C.~A. Weibel.
\newblock Mayer-{V}ietoris sequences and module structures on {$NK\sb\ast $}.
\newblock In {\em Algebraic $K$-theory, Evanston 1980 (Proc. Conf.,
  Northwestern Univ., Evanston, Ill., 1980)}, volume 854 of {\em Lecture Notes
  in Math.}, pages 466--493. Springer, Berlin, 1981.

\bibitem{Weibel(1989)}
C.~A. Weibel.
\newblock Homotopy algebraic {$K$}-theory.
\newblock In {\em Algebraic $K$-theory and algebraic number theory (Honolulu,
  HI, 1987)}, volume~83 of {\em Contemp. Math.}, pages 461--488. Amer. Math.
  Soc., Providence, RI, 1989.

\end{thebibliography}

\end{document}